\documentclass[12pt,a4paper]{article}
\usepackage[all,poly,color]{xy}
\usepackage{amsmath,amssymb, amsthm, amsfonts}
\usepackage[pagebackref,colorlinks]{hyperref}
\usepackage{enumerate}
\usepackage{comment}
     
\theoremstyle{plain}
\newtheorem {lemma}{Lemma}
\newtheorem {proposition}[lemma]{Proposition}
\newtheorem {theorem}[lemma]{Theorem}

\theoremstyle{definition}
\newtheorem {definition}[lemma]{Definition}
\newtheorem {remark}[lemma]{Remark}
\newtheorem {example}[lemma]{Example}

\newcommand{\N}{\mathbb{N}}

\newcommand{\reg}{\operatorname{reg}}

\newcommand{\End}{\operatorname{End}}
\newcommand{\sink}{\operatorname{sink}}

\newcommand{\op}{\operatorname{op}}
\newcommand{\Es}{\mathcal{S}}

\newcommand{\id}{\operatorname{id}}
\newcommand{\ghost}{\operatorname{ghost}}

\DeclareMathOperator{\RG}{\mathbf{RG}}
\DeclareMathOperator{\ABS}{\mathbf{ABS}}
\DeclareMathOperator{\ERG}{\mathbf{ERG}}
\DeclareMathOperator{\EABS}{\mathbf{EABS}}
\DeclareMathOperator{\MOD}{\mathbf{Mod}}

\title{Modules for Leavitt path algebras via extended algebraic branching systems}

\author{Raimund Preusser}
\AtEndDocument{\bigskip{\footnotesize%
  \textsc{Chebyshev Laboratory, St. Petersburg State University, Russia} \par  
  \textit{E-mail address:} \texttt{raimund.preusser@gmx.de} \par
}}
\date{}

\begin{document}
\maketitle

\begin{abstract} 
For a graph $E$, we introduce the notion of an extended $E$-algebraic branching system, generalising the notion of an $E$-algebraic branching system introduced by Gon\c{c}alves and Royer. We classify the extended $E$-algebraic branching systems and show that they induce modules for the corresponding Leavitt path algebra $L(E)$. Among these modules we find a class of nonsimple modules whose endomorphism rings are fields.  
\end{abstract}
\let\thefootnote\relax
\let\thefootnote\relax

\section{Introduction}
Leavitt path algebras are associative algebras constructed from directed graphs. They were introduced by G. Abrams and G. Aranda Pino in 2005 \cite{aap05} and independently by P. Ara, M. Moreno and E. Pardo in 2007 \cite{AMP}. 
The field of Leavitt path algebras has connections to many other branches of mathematics, such as functional analysis, symbolic dynamics, K-theory and noncommutative geometry, cf. \cite{AASbook}.

There have been a substantial number of papers devoted to the representation theory of Leavitt path algebras. P. Ara and M. Brustenga proved that the category of modules for a Leavitt path algebra $L(E)$ of a finite graph $E$ is equivalent to a quotient category of the category of modules for the path algebra $P(E)$ \cite{AB}. 
D. Gon\c{c}alves and D. Royer obtained modules for Leavitt path algebras by introducing the notion of an algebraic branching system \cite{GR}. X. Chen used infinite paths in $E$ to obtain simple modules for the Leavitt path algebra $L(E)$ \cite{C}. Numerous work followed, noteworthy the work of P. Ara and K. Rangaswamy producing new simple modules associated to infinite emitters and characterising those algebras which have countably (finitely) many distinct isomorphism classes of simple modules \cite{ARa,R-2,R-3}. G. Abrams, F. Mantese and A. Tonolo studied the projective resolutions for these simple modules \cite{AMT}, and P. \'Anh and T. Nam  provided another way to describe the so-called Chen and Rangaswamy simple modules \cite{nam}. In the recent paper \cite{HPS}, simple modules for \textit{weighted} Leavitt path algebras were found. Moreover, the algebraic branching systems for unweighted graphs and the modules defined by them were classified.

In this paper, we obtain modules for a Leavitt path algebra $L(E)$ over a field $K$ by introducing the notion of an \textit{extended} $E$-algebraic branching system. As a motivating example consider the graph
\[E:\xymatrix{ v\ar@[yellow]@(dr,ur)_{e}\ar@[green]@(dl,ul)^{d}}.\]
Set
\begin{align*}
X:=\{v\}\cup \{p,p^*\mid p\in E^{\geq 1}\}\cup\{pq^*\mid ~&p=p_1\dots p_m, q=q_1\dots q_n\in E^{\geq 1},(p_m,q_n)\neq (d,d)\}.
\end{align*}
Then $X$ is a linear basis for $L(E)$ by \cite[Theorem 1]{zelmanov1}. For any $i\in \{d,e,d^*,e^*\}$ let $Y_i$ be the subset of $X$ consisting of all elements whose last letter is $i$. Define the maps 
\begin{align*}
\rho_d:&~X\setminus (Y_{d^*}\cup Y_{e^*})\to Y_d,~x\mapsto xd;\\
\rho_e:&~X\setminus (Y_{d^*}\cup Y_{e^*})\to Y_e,~x\mapsto xe;\\
\rho_{d^*}:&~X\setminus Y_{d}\to Y_{d^*},~x\mapsto xd^*;\\
\rho_{e^*}:&~X\to Y_{e^*},~x\mapsto xe^{*}.
\end{align*}
Here we use the convention that $vi=i$ for any $i\in \{d,e,d^*,e^*\}$. Clearly the maps $\rho_d$, $\rho_e$, $\rho_{d^*}$ and $\rho_{e^*}$ are bijections. Let $\pi:L(E)\to \End(L(E))^{\op}$ be the right regular representation of $L(E)$. Then for any $x\in X$ one has
\begin{align}
\pi(v)(x)&=x,
\\
\pi(d)(x)&=\begin{cases}
\rho_{d}(x),&\text{ if }x\in X\setminus(Y_{d^*}\cup Y_{e^*}),\\
\rho_{d^*}^{-1}(x),&\text{ if }x\in Y_{d^*},\\
0,&\text{ if }x\in Y_{e^*},
\end{cases}
\\
\pi(e)(x)&=\begin{cases}
\rho_{e}(x),&\text{ if }x\in X\setminus(Y_{d^*}\cup Y_{e^*}),\\
0,&\text{ if }x\in Y_{d^*},\\
\rho_{e^*}^{-1}(x),&\text{ if }x\in Y_{e^*},
\end{cases}
\\
\pi(d^*)(x)&=\begin{cases}
\rho_{d^*}(x),&\text{ if }x\in X\setminus Y_{d},\\
\rho^{-1}_{d}(x)-\rho_{e^*}(\rho_{e}(\rho_{d}^{-1}(x))),&\text{ if }x\in Y_{d},
\end{cases}
\\
\pi(e^*)(x)&=\rho_{e^*}(x).
\end{align}
More generally, suppose that
\begin{itemize}
\item $Y_d$, $Y_e$, $Y_{d^*}$, $Y_{e^*}$ are \textit{any} pairwise disjoint subsets of a set $X$ and 
\item $\rho_d:X\setminus (Y_{d^*}\cup Y_{e^*})\to Y_d$, $\rho_e:X\setminus (Y_{d^*}\cup Y_{e^*})\to Y_e$, $\rho_{d^*}:X\setminus Y_{d}\to Y_{d^*}$, $\rho_{e^*}:X\to Y_{e^*}$ are \textit{any} bijections.
\end{itemize}
Let $V$ denote the $K$-vector space with basis $X$. Then it follows from the universal property of $L(E)$ that there is a unique representation $\pi:L(E)\to \End(V)^{\op}$ such that (1)-(5) hold. We call $\Es=(X,\{Y_i\}_{i\in \{d,e,d^*,e^*\}},\{\rho_j\}_{j\in \{d,e,d^*,e^*\}})$ an \textit{extended $E$-algebraic branching system} and the corresponding representation $\pi:L(E)\to \End(V)^{\op}$ the {\it representation of $L(E)$ associated to $\Es$}. The usual $E$-algebraic branching systems are in 1-1 correspondence to the extended $E$-algebraic branching systems $\Es=(X,\{Y_i\},\{\rho_j\})$ which have the property that $Y_d=Y_e=\emptyset$ (for details see Example \ref{exusual}). The notion of an extended $E$-algebraic branching system extends naturally to arbitary row-finite graphs $E$, see Section 4. 

The rest of the paper is organised as follows. In Section 2, we recall some graph-theoretic notions and the definition of a Leavitt path algebra. In Section 3, we recall the notion of an algebraic branching system and the notion of a representation graph (the latter was recently introduced in \cite{HPS}). Moreover, we recall the result from \cite[A.5]{HPS} that for a graph $E$, the categories $\ABS(E)$ of $E$-algebraic branching systems and $\RG(E)$ of representation graphs for $E$ are equivalent. In Section 4, we introduce the main notion of this paper, namely the notion of an extended algebraic branching system. Moreover, we introduce the notion of an extended representation graph and show that for a given graph $E$,  the categories $\EABS(E)$ of extended $E$-algebraic branching systems and $\ERG(E)$ of extended representation graphs for $E$ are equivalent. In Section 5, we classify the extended algebraic branching systems. For the usual algebraic branching systems we recover the classification obtained in \cite[Section 4]{HPS}. Suppose that $M$ is a module of an algebra over a field. Recall that Schur's lemma says that if $M$ is simple, then the endomorphism ring of $M$ is a skew field. In Section 6, we obtain for any Leavitt path algebra $L(E)$ over a field $K$ a class of nonsimple $L(E)$-modules whose endomorphism rings are isomorphic to $K$ (hence these modules constitute counter-examples to the converse of Schur's lemma).


\section{Preliminaries}
Throughout the paper $K$ denotes a fixed field. Rings and algebras are associative but not necessarily commutative or unital. 
$\N$ denotes the set of positive integers. 

\subsection{Graphs}

A {\it (directed) graph} is a quadruple $E=(E^0,E^1,s,r)$ where $E^0$ and $E^1$ are sets and $s,r:E^1\rightarrow E^0$ maps. The elements of $E^0$ are called {\it vertices} and the elements of $E^1$ {\it edges}. If $e$ is an edge, then $s(e)$ is called its {\it source} and $r(e)$ its {\it range}. If $v$ is a vertex and $e$ an edge, we say that $v$ {\it emits} $e$ if $s(e)=v$, and $v$ {\it receives} $e$ if $r(e)=v$. A vertex is called a {\it source} if it receives no edges, a {\it sink} if it emits no edges, and an \textit{infinite emitter} if it emits infinitely many edges. A vertex is called {\it regular} if it is neither a sink nor an infinite emitter. The subset of $E^0$ consisting of all sinks is denoted by $E^0_{\sink}$, and the subset consisting of all regular vertices by $E^0_{\reg}$. A graph is called \textit{nonempty} if it contains at least one vertex, and {\it row-finite} if it does not contain any infinite emitters. 

Let $E$ and $F$ be graphs. A {\it graph homomorphism} $\phi: E\to F$ consists of two maps $\phi^0 : E^0\to F^0$ and $\phi^1 : E^1\to F^1$ such that $s(\phi^1(e)) = \phi^0(s(e))$ and $r(\phi^1(e)) = \phi^0(r(e))$ for any $e\in E^1$. If $v\in E^0$ and $e\in E^1$, we will usually write $\phi(v)$ instead of $\phi^0(v)$ and $\phi(e)$ instead of $\phi^1(e)$. A graph $G$ is called a {\it subgraph} of a graph $E$ if $G^0\subseteq E^0$, $G^1\subseteq E^1$, $s_G=s_E|_{G^0}$ and $r_G=r_E|_{G^0}$.

Let $E=(E^0,E^1,s_E,r_E)$ and $F=(F^0,F^1,s_F,r_F)$ denote graphs. The graph $E\sqcup F=((E\sqcup F)^0,(E\sqcup F)^1,s,r)$ where $(E\sqcup F)^0=E^0\sqcup F^0$, $(E\sqcup F)^1=E^1\sqcup F^1$, $s|_{E^0}=s_E$, $s|_{F^0}=s_F$, $r|_{E^0}=r_E$ and $r|_{F^0}=r_F$ is called the \textit{disjoint union} of $E$ and $F$.

Let $E$ be a graph. The graph $E_d=(E_d^0, E_d^1, s_d, r_d)$ where $E_d^0=E^0$, $E_d^1=\{e,e^*\mid e\in E^1\}$, and $s_d(e)=s(e),~r_d(e)=r(e),~s_d(e^*)=r(e),~r_d(e^*)=s(e)$ for any $e\in E^1$ is called the {\it double graph} of $E$. We sometimes refer to the edges in the graph $E$ as {\it real edges} and to the additional edges in $E_d$ as {\it ghost edges}. The graph $\overline{E}$ obtained from $E_d$ by removing the real edges is called the \textit{inverse graph} of $E$.

A {\it (finite) path} in a graph $E$ is a finite, nonempty word $p=x_1\dots x_n$ over the alphabet $E^0\cup E^1$ such that either $x_i\in E^1~(i=1,\dots,n)$ and $r(x_i)=s(x_{i+1})~(i=1,\dots,n-1)$ or $n=1$ and $x_1\in E^0$. By definition, the {\it length} $|p|$ of $p$ is $n$ in the first case and $0$ in the latter case. The set of all paths in $E$ of length $n$ (resp. $\geq n$) is denoted by $E^n$ (resp. $E^{\geq n}$). 
For a path $p=x_1\dots x_n\in E^n$ we set $s(p):=s(x_1)$ and $r(p):=r(x_n)$ using the convention $s(v)=r(v)=v$ for any $v\in E^0$.

A {\it closed path} is a path $p\in E^{\geq 1}$ such that $s(p)=r(p)$. 
A {\it cycle} is a closed path $x_1\dots x_n$ such that $s(x_i)\neq s(x_j)$ for any $i\neq j$. An edge $e\in E^1$ is called an {\it exit} of a cycle $x_1\dots x_n$ if there is an $i\in \{1,\dots,n\}$ such that $s(e)=s(x_i)$ and $e\neq x_i$. A graph is called \textit{cyclic} if it contains a cycle, and \textit{acyclic} otherwise. 

A {\it (left-)infinite path} in a graph $E$ is a left-infinite word $p=\dots x_3x_2x_1$ over the alphabet $E^0\cup E^1$ such that for any $n\in\N$ the suffix $x_n\dots x_1$ is a path in $E$. We set $|p|:=\infty$ and $r(p):=r(x_1)$. The set of infinite paths in $E$ is denoted by $E^\infty$.

We say that a graph $E$ is \textit{connected} if for any $u,v\in E^0$ there is a path $p$ in the double graph $E_d$ such that $s_d(p)=u$ and $r_d(p)=v$. A maximal connected subgraph of $E$ is called a \textit{connected component} of $E$. Clearly any graph is the disjoint union of its connected components.

\subsection{Leavitt path algebras}
\begin{definition}\label{deflpa}
Let $E$ be a graph. The $K$-algebra $L(E)$ presented by the generating set $\{v,e,e^*\mid v\in E^0,e\in E^1\}$ and the relations
\begin{enumerate}[(i)]
\item $uv=\delta_{uv}u\quad(u,v\in E^0)$,
\item $s(e)e=e=er(e),~r(e)e^*=e^*=e^*s(e)\quad(e\in E^1)$,
\item $e^*f= \delta_{ef}r(e)\quad(e,f\in E^1)$,
\item $\sum_{e\in s^{-1}(v)}ee^*= v\quad(v\in E_{\reg}^0)$
\end{enumerate}
is called the {\it Leavitt path algebra} of $E$. 
\end{definition}




Let $E$ be a graph and $A$ a $K$-algebra. An {\it $E$-family} in $A$ is a subset $\{\sigma_v,\sigma_e,\sigma_{e^*}\mid v\in E^0, e\in E^1\}\subseteq A$ such that
\begin{enumerate}[(i)]
\item $\sigma_u\sigma_v=\delta_{uv}\sigma_u\quad(u,v\in E^0)$, 
\item
$\sigma_{s(e)}\sigma_{e}=\sigma_{e}=\sigma_{e}\sigma_{r(e)},~\sigma_{r(e)}\sigma_{e^*}=\sigma_{e^*}=\sigma_{e^*}\sigma_{s(e)}\quad(e\in E^1)$,
\item $\sigma_{e^*}\sigma_{f}= \delta_{ef}\sigma_{r(e)}\quad(e,f\in E^1)$ and
\item $\sum_{e\in s^{-1}(v)}\sigma_{e}\sigma_{e^*}= \sigma_{v}\quad(v\in E^0_{\reg})$.
\end{enumerate}
By the defining relations of $L(E)$, there exists a unique $K$-algebra homomorphism $\pi: L(E)\rightarrow A$ such that $\pi(v)=\sigma_v$, $\pi(e)=\sigma_{e}$ and $\pi(e^*)=\sigma_{e^*}$ for all $v\in E^0$ and $e\in E^1$. We will refer to this as the {\it universal property of $L(E)$}.

\section{Algebraic branching systems}
Until the end of this article $E$ denotes a fixed row-finite graph. For any $v\in E^0_{\reg}$ we choose an edge $e^v\in s^{-1}(v)$. The edges $e^v~(v\in E^0_{\reg})$ are called {\it special} and the other edges in $E$ are called {\it nonspecial}.

\subsection{Algebraic branching systems}

\begin{definition}\label{defabs}
Let
\begin{itemize}
\item $\{X_v\}_{v\in E^0}$ be a family of pairwise disjoint subsets of a set $X$,
\item $\{Y_{e^*}\}_{e\in E^1}$ a family of pairwise disjoint sets such that $X_v=\bigcup_{e\in s^{-1}(v)}Y_{e^*}$ for any $v\in E_{\reg}^0$ and 
\item $\{\rho_{e^*}\}_{e\in E^1}$ a family of maps such that $\rho_{e^*}:X_{r(e)} \to Y_{e^*}$ is a bijection for any $e\in E^1$.
\end{itemize}
Then $(X,\{X_v\}_{v\in E^0},\{Y_{e^*}\}_{e\in E^1},\{\rho_{e^*}\}_{e\in E^1})$
is called an {\it $E$-algebraic branching system}. 
In this paper all $E$-algebraic branching systems are assumed to be \textit{saturated}, i.e. $X=\bigcup_{v\in E^0}X_v$.
\end{definition}

\begin{remark}
We use the notation $Y_{e^*}$ (respectively $\rho_{e^*}$) instead of $Y_{e}$ (respectively $\rho_{e}$) because in this way it is easier to see that algebraic branching systems correspond to a certain type of extended algebraic branching systems, cf. Section 4.
\end{remark}

We denote by $\ABS(E)$ the category of $E$-algebraic branching systems. A morphism $\alpha:\Es\to\Es'$ between two objects $\Es=(X,\{X_v\},\{Y_{e^*}\},\{\rho_{e^*}\})$ and $\Es'=(X',\{X'_v\},\{Y'_{e^*}\},$ $\{\rho'_{e^*}\})$ in $\ABS(E)$ is a map $\alpha: X\to X'$ such that $\alpha(X_v)\subseteq X'_v$ and $\alpha (Y_{e^*})\subseteq Y'_{e^*}$ for any $v\in E^0$ and $e\in E^1$, and $\alpha$ is compatible with the bijections inside $\Es$ and $\Es'$. 

Suppose that $\Es=(X,\{X_v\},\{Y_{e^*}\},\{\rho_{e^*}\})$ is an $E$-algebraic branching system. Let $V=V(\Es)$ denote the $K$-vector space with basis $X$. For any $v\in E^0$ and $e\in E^1$ define endomorphisms $\sigma_v,\sigma_e,\sigma_{e^*}\in \End(V)$ by 

\begin{align*}
\sigma_v(x)&=\begin{cases}
x,\quad\quad\hspace{0.1cm}&\text{ if }x\in X_v,\\
0,&\text{ otherwise,}
\end{cases}\\
\sigma_{e}(x)&
=\begin{cases}
\rho^{-1}_{e^*}(x),&\text{ if }x\in Y_{e^*},\\
0,&\text{ otherwise,}
\end{cases}\\
\sigma_{e^*}(x)&=\begin{cases}
\rho_{e^*}(x),\hspace{0.1cm}&\text{ if }x\in X_{r(e)},\\
0,&\text{ otherwise,}
\end{cases}
\end{align*}
for any $x\in X$. It follows from the universal property of $L(E)$ that there is a $K$-algebra homomorphism $\pi:L(E)\to \End(V)^{\op}$ such that $\pi(v)=\sigma_v$, $\pi(e)=\sigma_e$ and $\pi(e^*)=\sigma_{e^*}$ for any $v\in E^0$ and $e\in E^1$. We call $\pi:L(E)\to \End(V)^{\op}$ the {\it representation of $L(E)$ associated to $\Es$}. 

\begin{remark}
Let $\Es=(X,\{X_v\},\{Y_{e^*}\},\{\rho_{e^*}\})$ be an $E$-algebraic branching system. Let $V$ denote the $K$-vector space with basis $X$ and $V'$ the $K$-vector space consisting of all maps from $X$ to $K$. In \cite{GR}, the representation of $L(E)$ associated to $\Es$ was defined using $V'$ instead of $V$, but the possibility of using $V$ was mentioned in \cite[Remark 2.3]{GR}. In \cite{C} the representation of $L(E)$ associated to $\Es$ was defined using $V$. 
\end{remark}

If $\Es$ is an object in $\ABS(E)$, then the vector space $V(\Es)$ becomes a right $L(E)$-module by defining $a.b:=\pi(b)(a)$ for any $a\in V(\Es)$ and $b\in L(E)$. We call this module the \textit{$L(E)$-module associated to $\Es$}. If $\alpha:\Es\to\Es'$ is a morphism in $\ABS(E)$, let $V(\alpha):V(\Es)\to V(\Es')$ be the module homomorphism such that $V(\alpha)(x)=\alpha(x)$ for any $x\in X$. We obtain a functor
\[V:\ABS(E)\to \MOD(L(E))\]
where $\MOD(L(E))$ denotes the category of right $L(E)$-modules.

\subsection{Representation graphs}
One can visualise branching systems using representation graphs, which are defined below. Recall from \S 2.1 that $\overline{E}$ denotes the inverse graph of $E$.

\begin{definition}\label{defrg}
A {\it representation graph} for $E$ is a pair $(F,\phi)$, where $F$ is a graph and $\phi:F\rightarrow \overline{E}$ is a graph homomorphism such that the following hold for any $w\in F^0$.
\begin{enumerate}[(i)]
\item If $\phi(w)$ receives an edge in $\overline{E}$, then $|r^{-1}(w)|=1$.
\item $\phi$ maps $s^{-1}(w)$ 1-1 onto $s_{\overline{E}}^{-1}(\phi(w))$. 
\end{enumerate}
\end{definition}

\begin{example}\label{ex-1}
Suppose that $E$ is the graph $\xymatrix{ v\ar@[yellow]@(dr,ur)_{e}\ar@[green]@(dl,ul)^{d}}$ whose inverse graph $\overline{E}$ is $\xymatrix{ v\ar@[blue]@(dr,ur)_{e^*}\ar@[red]@(dl,ul)^{d^*}}$. Then a representation graph $(F,\phi)$ for $E$ is given by
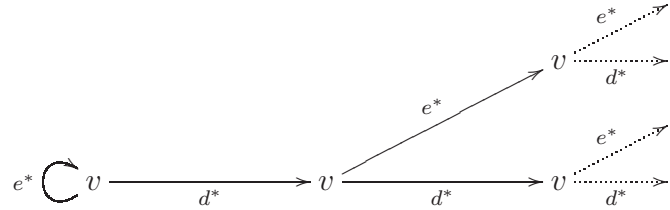
\begin{equation*}
\xymatrix@R=0.5cm@C=0.5cm{
&&&&&&&&&&&\ar@[blue]@{<.}[dll]_{e^*}\\
&&&&&&&&&v\ar@[blue]@{<-}[ddllll]_{e^*}&&\ar@[red]@{<.}[ll]^{d^*}\\
&&&&&&&&&&&\ar@[blue]@{<.}[dll]_{e^*}\\
&v\ar@[blue]@(dl,ul)^{e^*}&&&&v\ar@[red]@{<-}[llll]^{d^*}&&&&v\ar@[red]@{<-}[llll]^{d^*}&&\ar@[red]@{<.}[ll]^{d^*}}
\end{equation*}
where a label of a vertex (resp. edge) indicates the image of that vertex (resp. edge) under $\phi$.
\end{example}

We denote by $\RG(E)$ the category of representation graphs for $E$. A morphism $\alpha:(F,\phi)\to (G,\psi)$ in $\RG(E)$ is a graph homomorphism $\alpha:F\to G$ such that $\psi\circ \alpha=\phi$.

Suppose that $(F,\phi)$ is an object in $\RG(E)$. Let $W=W(F,\phi)$ be the $K$-vector space with basis $F^0$. For any $v\in E^0$ and $e\in E^1$ define endomorphisms $\sigma_v,\sigma_e,\sigma_{e^*}\in \End(W)$ by
\begin{align*}
\sigma_v(w)&=\begin{cases}w,\quad\quad\hspace{0.05cm}&\text{if }\phi(w)=v,\\0,& \text{otherwise}, \end{cases}\\
\sigma_e(w)&=\begin{cases}s(f),\quad&\text{if }\exists f\in r^{-1}(w):\phi(f)=e^*,\\0,& \text{otherwise}, \end{cases}\\
\sigma_{e^*}(w)&=\begin{cases}r(f),\quad&\text{if }\exists f\in s^{-1}(w):\phi(f)=e^*,\\0,& \text{otherwise},
 \end{cases}
\end{align*}
for any $w\in F^0$. It follows from the universal property of $L(E)$ that there is a $K$-algebra homomorphism $\pi:L(E)\to \End(W)^{\op}$ such that $\pi(v)=\sigma_v$, $\pi(e)=\sigma_e$ and $\pi(e^*)=\sigma_{e^*}$ for any $v\in E^0$ and $e\in E^1$. We call $\pi:L(E)\to \End(W)^{\op}$ the {\it representation of $L(E)$ associated to $(F,\Phi)$}. 

If $(F,\phi)$ is an object in $\RG(E)$, then the vector space $W(F,\phi)$ becomes a right $L(E)$-module by defining $a.b:=\pi(b)(a)$ for any $a\in W(F,\phi)$ and $b\in L(E)$. We call this module the \textit{$L(E)$-module associated to $(F,\phi)$}. If $\alpha:(F,\phi)\to (G,\psi)$ is a morphism in $\RG(E)$, let $W(\alpha):W(F,\phi)\to W(G,\psi)$ be the module homomorphism such that $W(\alpha)(w)=\alpha(w)$ for any $w\in F^0$. We obtain a functor
\[W:\RG(E)\to \MOD(L(E)).\]

\begin{example}\label{ex0}
Suppose again that $E$ is the graph $\xymatrix{ v\ar@[yellow]@(dr,ur)_{e}\ar@[green]@(dl,ul)^{d}}$ whose inverse graph $\overline{E}$ is $\xymatrix{ v\ar@[blue]@(dr,ur)_{e^*}\ar@[red]@(dl,ul)^{d^*}}$, and that $(F,\phi)$ is the representation graph for $E$ given by
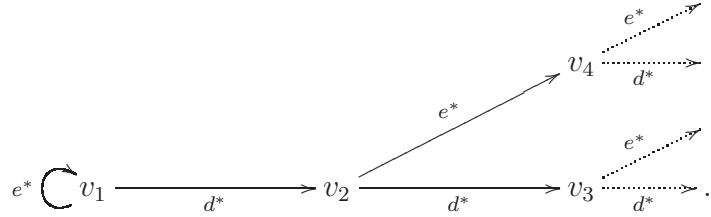
\begin{equation*}
\xymatrix@R=0.5cm@C=0.5cm{
&&&&&&&&&&&\ar@[blue]@{<.}[dll]_{e^*}\\
&&&&&&&&&v_4\ar@[blue]@{<-}[ddllll]_{e^*}&&\ar@[red]@{<.}[ll]^{d^*}\\
&&&&&&&&&&&\ar@[blue]@{<.}[dll]_{e^*}\\
&v_1\ar@[blue]@(dl,ul)^{e^*}&&&&v_2\ar@[red]@{<-}[llll]^{d^*}&&&&v_3\ar@[red]@{<-}[llll]^{d^*}&&\ar@[red]@{<.}[ll]^{d^*}.}
\end{equation*}
Recall that $W(F,\Phi)$ is the $K$-vector space with basis $F^0$. The action of $L(E)$ on $W(F,\Phi)$ ``slides'' the vertices of $F$ along the edges. For example, $v_1. d^*=v_2$, $v_2.d^*=v_3$, $v_2. d=v_1$ and $v_2. e=0$. The $L(E)$-module $W(F,\Phi)$ is isomorphic to the Chen simple module defined by the right-infinite path $eee\dots $, see \cite[Section 4]{HPS}. 
\end{example}

\subsection{Algebraic branching systems vs. representation graphs}
To any object $\Es=(X,\{X_v\},\{Y_{e^*}\},\{\rho_{e^*}\})$ in $\ABS (E)$ we associate the object $\eta(\Es)=(F,\phi)$ in $\RG(E)$ defined by
\begin{itemize}
\item $F^0=X$,
\item $F^1=\{f_x\mid x\in \bigcup_{e\in E^1}Y_{e^*}\}$,
\item $s_F(f_x)=\rho_{e^*}^{-1}(x)$ if $x\in Y_{e^*}$,
\item $r_F(f_x)=x$,
\item $\phi^0(x)=v$ if $x\in X_v$,
\item $\phi^1(f_x)=e^*$ if $x\in Y_{e^*}$.
\end{itemize}
To any morphism $\alpha:\Es=(X,\{X_v\},\{Y_{e^*}\},\{\rho_{e^*}\})\to \Es'=(X',\{X'_v\},\{Y'_{e^*}\},\{\rho'_{e^*}\})$ in $\ABS (E)$ we associate the morphism $ \eta(\alpha):\eta(\Es)\to \eta(\Es')$ in $\RG(E)$ defined by $ \eta(\alpha)^0(x)=\alpha(x)$ for any $x\in X$ and $ \eta(\alpha)^1(f_x)=f'_{\alpha(x)}$ for any $x\in \bigcup_{e\in E^1}Y_{e^*}$. In this way we obtain a functor $\eta:\ABS (E)\to \RG(E)$. 

Conversely, to any object $(F,\phi)$ in $\RG(E)$ we associate an object $ \theta(F,\phi)=(X,\{X_v\},\{Y_{e^*}\},\{\rho_{e^*}\})$ in $\ABS (E)$ defined by
\begin{itemize}
\item $X=F^0$,
\item $X_v=\{w\in F^0\mid \phi(w)=v\}$ for any $v\in E^0$,
\item $Y_{e^*}=\{w\in F^0\mid \exists f\in r_F^{-1}(w):\phi(f)=e^*\}$ for any $e\in E^1$,
\item $\rho_{e^*}(w)=r(f)$ for any $e\in E^1$ and $w\in X_{r(e)}$, where $f$ is the unique edge in $s^{-1}(w)\cap \phi^{-1}(e^*)$. 
\end{itemize}
To any morphism $\alpha:(F,\phi)\to (G,\psi)$ in $\RG(E)$ we associate the morphism $ \theta(\alpha): \theta(F,\phi)\to \theta(G,\psi)$ in $\ABS(E)$ defined by $ \theta(\alpha)(x)=\alpha(x)$ for any $x\in X=F^0$. In this way we obtain a functor $ \theta:\RG(E)\to \ABS (E)$. 

We leave it to the reader to check that $ \theta\circ\eta=\id_{\ABS (E)}$ and $ \eta\circ   \theta\cong\id_{\RG(E)}$. Hence the categories $\RG(E)$ and $\ABS (E)$ are equivalent. Moreover, the diagrams
\[\xymatrix{\ABS (E)\ar[dr]^V&\\&\MOD L(E)& \text{and}\\\RG(E)\ar[uu]^{\theta}\ar[ur]_W&\\}
\quad\quad
\xymatrix{\ABS (E)\ar[dd]_{\eta}\ar[dr]^V&\\&\MOD L(E)\\\RG(E)\ar[ur]_W&\\}
\]
are commutative.

\section{Extended algebraic branching systems}
\subsection{Extended algebraic branching systems}
\begin{definition}\label{defeabs}
Let 
\begin{itemize}
\item $\{X_v\}_{v\in E^0}$ be a family of pairwise disjoint subsets of a set $X$,
\item $\{Y_i\}_{i\in E^1\cup (E^1)^*}$ a family of pairwise disjoint sets such that
$Y_e\subseteq X_{r(e)}$ and $Y_{e^*}\subseteq X_{s(e)}$ for any $e\in E^1$,
\item $\{\rho_j\}_{j\in E^1\cup (E^1)^*}$ a family of maps such that 
\begin{enumerate}[(i)]
\item $\rho_{e}: X_{s(e)}\setminus\bigcup_{f\in s^{-1}(s(e))}Y_{f^*}\to Y_{e}$ is a bijection for any $e\in E^1$,
\item $\rho_{e^*}: X_{r(e)}\to Y_{e^*}$ is a bijection for any nonspecial $e\in E^1$,
\item $\rho_{e^*}: X_{r(e)}\setminus Y_{e}\to Y_{e^*}$ is a bijection for any special $e\in E^1$.
\end{enumerate}
\end{itemize}
Then \[(X, \{X_v\}_{v\in E^0}, \{Y_i\}_{i\in E^1\cup (E^1)^*},\{\rho_j\}_{j\in E^1\cup (E^1)^*})\]
is called an {\it extended $E$-algebraic branching system}. In this paper all extended $E$-algebraic branching systems are assumed to be \textit{saturated}, i.e. $X=\bigcup_{v\in E^0}X_v$.
\end{definition}

\begin{example}\label{exusual}
Suppose that $(X, \{X_v\}_{v\in E^0}, \{Y_{e^*}\}_{e\in  E^1},\{\rho_{e^*}\}_{e\in E^1})$ is a usual $E$-algebraic branching system. For any $e\in E^1$ let $Y_e=\emptyset$ and $\rho_{e}: \emptyset=X_{s(e)}\setminus\bigcup_{f\in s^{-1}(s(e))}Y_{f^*}\to Y_{e}=\emptyset$ be the unique map. Then $(X, \{X_v\}_{v\in E^0}, \{Y_i\}_{i\in  E^1\cup (E^1)^*},\{\rho_j\}_{j\in E^1\cup (E^1)^*})$ is an extended $E$-algebraic branching system. In this way one obtains a 1-1 correspondence between the usual $E$-algebraic branching systems and the extended $E$-algebraic branching systems which have the property that $Y_e=\emptyset$ for any $e\in E^1$.
\end{example}

\begin{example}
Set
\begin{align*}
X:=\{v\mid v\in E^0\}\cup\{p,p^*\mid p\in E^{\geq 1}\}\cup\{pq^*\mid \hspace{0.12cm}&p=e_1\dots e_m, q=f_1\dots f_n\in E^{\geq 1},r(p)=r(q)\\&\text {and either }e_m\neq f_n\text{ or }e_m=f_n\text{ is nonspecial}\}.
\end{align*}
Then $X$ is a linear basis of $L(E)$ by \cite[Theorem 1]{zelmanov1}. For any $v\in E^0$ let  $X_v$ be the subset of $X$ consisting of all elements whose range is $v$. For any $i\in E^1\cup (E^1)^*$ let $Y_i$ be the subset of $X$ consisting of all elements whose last letter is $i$. For any $j\in E^1\cup (E^1)^*$ let $\rho_j$ be the map that adds the letter $j$ to the end of a path respectively replaces a letter from $E^0$ by $j$ (the domain and codomain of $\rho_j$ are as in Definition \ref{defeabs}). One checks easily that $(X,\{X_v\},\{Y_i\},\{\rho_j\})$ is an extended $E$-algebraic branching system. 
\end{example}

We denote by $\EABS(E)$ the category of extended $E$-algebraic branching systems. A morphism $\alpha:\Es\to\Es'$ between two objects $\Es=(X,\{X_v\},\{Y_{i}\},\{\rho_{j}\})$ and $\Es'=(X',\{X'_v\},\{Y'_{i}\},$ $\{\rho'_{j}\})$ in $\EABS(E)$ is a map $\alpha: X\to X'$ such that $\alpha (X_{v})\subseteq X'_{v}$ and $\alpha (Y_{i})\subseteq Y'_{i}$ for any $v\in E^0$ and $i\in E^1\cup (E^1)^*$, and $\alpha$ is compatible with the bijections inside $\Es$ and $\Es'$ (in particular we require that $\alpha(X_{s(e)}\setminus\bigcup_{f\in s^{-1}(s(e))}Y_{f^*})\subseteq X'_{s(e)}\setminus\bigcup_{f\in s^{-1}(s(e))}Y'_{f^*}$ for any $e\in E^1$ and $\alpha(X_{r(e)}\setminus Y_{e})\subseteq X'_{r(e)}\setminus Y'_{e}$ for any special $e\in E^1$). 

Suppose that $\Es=(X,\{X_v\},\{Y_{i}\},\{\rho_{j}\})$ is an extended $E$-algebraic branching system. Let $V=V(\Es)$ denote the $K$-vector space with basis $X$. For any $v\in E^0$ and $e\in E^1$ define endomorphisms $\sigma_v,\sigma_e,\sigma_{e^*}\in \End(V)$ by 

\begin{align*}
\sigma_v(x)&=\begin{cases}
x,&\hspace{6.9cm}\text{ if }x\in X_v,\\
0,&\hspace{6.9cm}\text{ otherwise,}
\end{cases}\\
\sigma_{e}(x)&=\begin{cases}
\rho_{e}(x),&\hspace{5.95cm}\text{ if }x\in X_{s(e)}\setminus\bigcup_{f\in s^{-1}(s(e))}Y_{f^*},\\
\rho_{e^*}^{-1}(x),&\hspace{5.95cm}\text{ if }x\in Y_{e^*},\\
0,&\hspace{5.95cm}\text{ otherwise,}
\end{cases}\\
\sigma_{e^*}(x)&=\begin{cases}
\rho_{e^*}(x),&\text{ if }x\in X_{r(e)}\setminus Y_{e},\\
\rho_{e^*}(x),&\text{ if }x\in Y_{e}\text{ and }e\text{ is nonspecial},\\
\rho^{-1}_{e}(x)-\sum_{d\in s^{-1}(s(e))\setminus\{e\}}\rho_{d^*}(\rho_{d}(\rho_{e}^{-1}(x))),&\text{ if }x\in Y_{e}\text{ and }e\text{ is special},\\
0,&\text{ otherwise,}
\end{cases}
\end{align*}
for any $x\in X$. It follows from the universal property of $L(E)$ that there is a $K$-algebra homomorphism $\pi:L(E)\to \End(V)^{\op}$ such that $\pi(v)=\sigma_v$, $\pi(e)=\sigma_e$ and $\pi(e^*)=\sigma_{e^*}$ for any $v\in E^0$ and $e\in E^1$. We call $\pi:L(E)\to \End(V)^{\op}$ the {\it representation of $L(E)$ associated to $\Es$}. 

If $\Es$ is an object in $\EABS(E)$, then the vector space $V(\Es)$ becomes a right $L(E)$-module by defining $a.b:=\pi(b)(a)$ for any $a\in V(\Es)$ and $b\in L(E)$. We call this module the \textit{$L(E)$-module associated to $\Es$}. If $\alpha:\Es\to\Es'$ is a morphism in $\EABS(E)$,  let $V(\alpha):V(\Es)\to V(\Es')$ be the module homomorphism such that $V(\alpha)(x)=\alpha(x)$ for any $x\in X$. We obtain a functor
\[V:\EABS(E)\to \MOD(L(E)).\]

\subsection{Extended representation graphs}

One can visualise extended branching systems using extended representation graphs, which are defined below. Recall from \S 2.1 that $E_d$ denotes the double graph of $E$.

\begin{definition}\label{deferg}
An {\it extended representation graph} for $E$ is a pair $(F,\phi)$, where $F$ is a graph and $\phi:F\rightarrow E_d$ a graph homomorphism such that the following hold for any $w\in F^0$.
\begin{enumerate}[(i)]
\item $w$ is either a source or receives a unique edge $f_w$.
\item If either $w$ is a source or $\phi(f_w)$ is a nonspecial real edge, then $\phi$ maps $s^{-1}(w)$ 1-1 onto $s_d^{-1}(\phi(w))$.
\item If $\phi(f_w)$ is a special real edge, then $\phi$ maps $s^{-1}(w)$ 1-1 onto $s_d^{-1}(\phi(w))\setminus\{(\phi(f_w))^*\}$.
\item If $\phi(f_w)$ is a ghost edge, then $\phi$ maps $s^{-1}(w)$ 1-1 onto $s_d^{-1}(\phi(w))\cap (E^1)^*$.
\end{enumerate}
\end{definition}
\vspace{0.1cm}

\begin{example}\label{extrivial}
The usual representation graphs for $E$ are precisely the extended representation graphs $(F,\phi)$ for $E$ having the property that $\phi(f)$ is a ghost edge for any $f\in F^1$.
\end{example}

\begin{example}\label{ex1}
Suppose that $E$ is the graph $\xymatrix{ v\ar@[yellow]@(dr,ur)_{e}\ar@[green]@(dl,ul)^{d}}$ whose double graph $E_d$ is
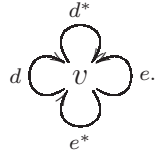
\begin{equation*}
\xymatrix{ v\ar@[yellow]@(dr,ur)_{e.}\ar@[green]@(dl,ul)^{d}\ar@[blue]@(dr,dl)^{e^*}\ar@[red]@(ul,ur)^{d^*}}
\end{equation*} 
As special edge choose $e=e^v$. Then an extended representation graph $(F,\phi)$ for $E$ is given by
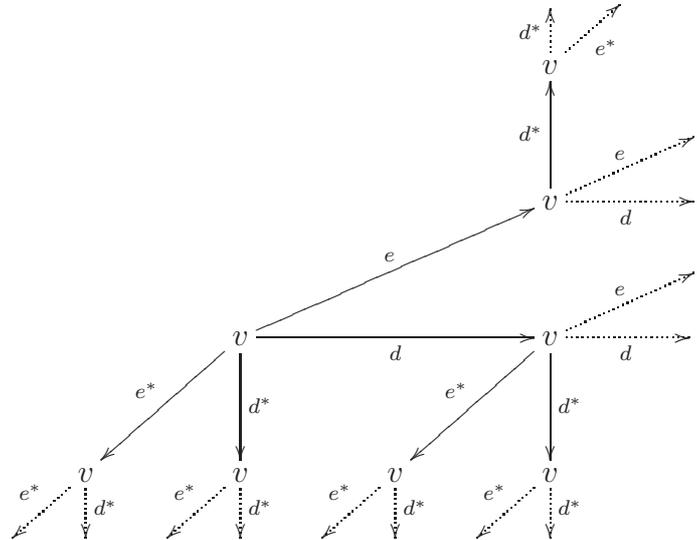
\begin{equation*}
\xymatrix@R=0.6cm@C=0.7cm{
&&&&&&&&&\\
&&&&&&&v\ar@[red]@{.>}[u]^{d^*}\ar@[blue]@{.>}[ur]_{e^*}&&\\
&&&&&&&&&\ar@[yellow]@{<.}[dll]_e\\
&&&&&&&v\ar@[yellow]@{<-}[ddllll]_e\ar@[red][uu]^{d^*}&&\ar@[green]@{<.}[ll]^d\\
&&&&&&&&&\ar@[yellow]@{<.}[dll]_e\\
&&&v\ar@[red][dd]^{d^*}\ar@[blue][ddll]_{e^*}&&&&v\ar@[green]@{<-}[llll]^d\ar@[red][dd]^{d^*}\ar@[blue][ddll]_{e^*}&&\ar@[green]@{<.}[ll]^d.\\
&&&&&&&&&\\
&v\ar@[red]@{.>}[d]^{d^*}\ar@[blue]@{.>}[dl]_{e^*}&&v\ar@[red]@{.>}[d]^{d^*}\ar@[blue]@{.>}[dl]_{e^*}&&v\ar@[red]@{.>}[d]^{d^*}\ar@[blue]@{.>}[dl]_{e^*}&&v\ar@[red]@{.>}[d]^{d^*}\ar@[blue]@{.>}[dl]_{e^*}&&\\
&&&&&&&&&\\
&&&&&&&&&}
\end{equation*}
\end{example}

We denote by $\ERG(E)$ the category of extended representation graphs for $E$. A morphism $\alpha:(F,\phi)\to (G,\psi)$ in $\ERG(E)$ is a graph homomorphism $\alpha:F\to G$ such that $\psi\circ \alpha=\phi$.

Suppose that $(F,\phi)$ is an object in $\ERG(E)$. Let $W=W(F,\phi)$ be the $K$-vector space with basis $F^0$. For any $v\in E^0$ and $e\in E^1$ define endomorphisms $\sigma_v,\sigma_e,\sigma_{e^*}\in \End(W)$ by
\begin{align*}
\sigma_v(w)&=\begin{cases}w,\quad\quad\quad\quad~&\text{if }\phi(w)=v,\\0,& \text{otherwise}, \end{cases}\\
\sigma_{e}(w)&=\begin{cases}
r(f),\quad&\quad~~~\text{ if }\exists f\in s^{-1}(w): \phi(f)=e,\\
s(f),\quad&\quad~~~\text{ if }\exists f\in r^{-1}(w): \phi(f)=e^*,\\
0,&\quad~~~\text{ otherwise,}
\end{cases}\\
\sigma_{e^*}(w)&=\begin{cases}
r(f),\quad&\text{ if }\exists f\in s^{-1}(w): \phi(f)=e^*,\\
s(f)-T,\quad&\text{ if }e\text{ is special and }\exists f\in r^{-1}(w): \phi(f)=e,\\
0,&\text{ otherwise,}
\end{cases}
\end{align*}
for any $w\in F^0$. Here $T=\sum_p r(p)$ where $p$ ranges over all paths in $F$ starting in $s(f)$ such that $\phi(p)=dd^*$ for some $d\in s^{-1}(s(e))\setminus\{e\}$. It follows from the universal property of $L(E)$ that there is a $K$-algebra homomorphism $\pi:L(E)\to \End(W)^{\op}$ such that $\pi(v)=\sigma_v$, $\pi(e)=\sigma_e$ and $\pi(e^*)=\sigma_{e^*}$ for any $v\in E^0$ and $e\in E^1$. We call $\pi:L(E)\to \End(W)^{\op}$ the {\it representation of $L(E)$ associated to $(F,\phi)$}. 

If $(F,\phi)$ is an object in $\ERG(E)$, then the vector space $W(F,\phi)$ becomes a right $L(E)$-module by defining $a.b:=\pi(b)(a)$ for any $a\in W(\Es)$ and $b\in L(E)$. We call this module the \textit{$L(E)$-module associated to $(F,\phi)$}. If $\alpha:(F,\phi)\to (G,\psi)$ is a morphism in $\ERG(E)$, let $W(\alpha):W(F,\phi)\to W(G,\psi)$ be the module homomorphism such that $W(\alpha)(w)=\alpha(w)$ for any $w\in F^0$. We obtain a functor
\[W:\ERG(E)\to \MOD(L(E)).\]

\subsection{Extended algebraic branching systems vs. extended representation graphs}
To any object $\Es=(X,\{X_v\},\{Y_{i}\},\{\rho_{j}\})$ in $\EABS (E)$ we associate the object $\eta(\Es)=(F,\phi)$ in $\ERG(E)$ defined by
\begin{itemize}
\item $F^0=X$,
\item $F^1=\{f_x\mid x\in \bigcup_{e\in E^1}(Y_e\cup Y_{e^*})\}$,
\item $s_F(f_x)=\rho_{e}^{-1}(x)$ if $x\in Y_{e}$, respectively $s_F(f_x)=\rho_{e^*}^{-1}(x)$ if $x\in Y_{e^*}$,
\item $r_F(f_x)=x$,
\item $\phi^0(x)=v$ if $x\in X_v$,
\item $\phi^1(f_x)=e$ if $x\in Y_{e}$, respectively $\phi^1(f_x)=e^*$ if $x\in Y_{e^*}$.
\end{itemize}
To any morphism $\alpha:\Es=(X,\{X_v\},\{Y_{i}\},\{\rho_{j}\})\to \Es'=(X',\{X'_v\},\{Y'_{i}\},\{\rho'_{j}\})$ in $\EABS (E)$ we associate the morphism $ \eta(\alpha):\eta(\Es)\to \eta(\Es')$ in $\ERG(E)$ defined by $ \eta(\alpha)^0(x)=\alpha(x)$ for any $x\in X$ and $ \eta(\alpha)^1(f_x)=f'_{\alpha(x)}$ for any $x\in \bigcup_{e\in E^1}(Y_e\cup Y_{e^*})$. In this way we obtain a functor $\eta:\EABS (E)\to \ERG(E)$. 

Conversely, to any object $(F,\phi)$ in $\ERG(E)$ we associate an object $ \theta(F,\phi)=(X,\{X_v\},\{Y_{i}\},\{\rho_{j}\})$ in $\EABS (E)$ defined by
\begin{itemize}
\item $X=F^0$,
\item $X_v=\{w\in F^0\mid \phi(w)=v\}$ for any $v\in E^0$,
\item $Y_{e}=\{w\in F^0\mid \exists f\in r_F^{-1}(w):\phi(f)=e\}$ for any $e\in E^1$,
\item $Y_{e^*}=\{w\in F^0\mid \exists f\in r_F^{-1}(w):\phi(f)=e^*\}$ for any $e\in E^1$,
\item $\rho_{e}(w)=r(f)$ for any $e\in E^1$, where $f$ is the unique edge in $s^{-1}(w)\cap \phi^{-1}(e)$,
\item $\rho_{e^*}(w)=r(f)$ for any $e\in E^1$, where $f$ is the unique edge in $s^{-1}(w)\cap \phi^{-1}(e^*)$. 
\end{itemize}
To any morphism $\alpha:(F,\phi)\to (G,\psi)$ in $\ERG(E)$ we associate the morphism $ \theta(\alpha): \theta(F,\phi)\to \theta(G,\psi)$ in $\EABS(E)$ defined by $ \theta(\alpha)(x)=\alpha(x)$ for any $x\in X=F^0$. In this way we obtain a functor $ \theta:\ERG(E)\to \EABS (E)$. 

We leave it to the reader to check that $ \theta\circ   \eta=\id_{\EABS (E)}$ and $ \eta\circ   \theta\cong\id_{\ERG(E)}$. Hence the categories $\ERG(E)$ and $\EABS (E)$ are equivalent. Moreover, the diagrams
\[\xymatrix{\EABS (E)\ar[dr]^V&\\&\MOD L(E)& \text{and}\\\ERG(E)\ar[uu]^{\theta}\ar[ur]_W&\\}
\quad\quad
\xymatrix{\EABS (E)\ar[dd]_{\eta}\ar[dr]^V&\\&\MOD L(E)\\\ERG(E)\ar[ur]_W&\\}
\]
are commutative.

\section{Classification of the extended algebraic branching systems}
The goal of this section is to classify the extended $E$-algebraic branching systems. In order to do so, it suffices to classify the extended representation graphs for $E$ (in view of \S4.3). 

If $(F,\phi)$ and $(G,\psi)$ are extended representation graphs for $E$, then their \textit{disjoint union} is the extended representation graph $(F\sqcup G, \phi\sqcup \psi)$ for $E$, where $(\phi\sqcup\psi)|_{F}=\phi$ and $(\phi\sqcup\psi)|_{G}=\psi$. If $(F,\phi)$ is an extended representation graph for $E$ and $\{F_i\}_{i\in\Lambda}$ are the connected components of $F$, then clearly $(F,\phi)$ is the disjoint union of the extended representation graphs $(F_i,\phi|_{F_i})$ for $E$. Hence it suffices to classify the \textit{connected} extended representation graphs for $E$.

\subsection{Basis paths}
\subsubsection{Finite basis paths}
Set
\begin{align*}
X:=\{v\mid v\in E^0\}\cup\{p,p^*\mid p\in E^{\geq 1}\}\cup\{pq^*\mid \hspace{0.12cm}&p=e_1\dots e_m, q=f_1\dots f_n\in E^{\geq 1},r(p)=r(q)\\&\text {and either }e_m\neq f_n\text{ or }e_m=f_n\text{ is nonspecial}\}.
\end{align*}
We consider the elements of $X$ as paths in $E_d$ and call them \textit{(finite) basis paths}. 
For any $v\in E^0$ we set $X_v:=\{x\in X\mid s_d(x)=v\}$. For a basis path $x=x_1\dots x_l\in X$ and $0\leq n \leq l$ we define $\tau_{\leq n}(x):=x_{1}\dots x_n\in X$ and $\tau_{>n}(x):=x_{n+1}\dots x_l\in X$. Here we use the convention $\tau_{\leq 0}(x)=s_d(x)$ and $\tau_{>l}(x)=r_d(x)$. 
 
For any graph $F$, $w\in F^0$ and $n\geq 0$ we set $F^{\geq n}_w:=\{p\in F^{\geq n} \mid s(p)=w\}$. Suppose $(F,\phi)$ is an extended representation graph for $E$. If $w\in F^0$ and $f\in F^1$, we say that \textit{$w$ lies over $\phi(w)$} and \textit{$f$ lies over $\phi(f)$}. If $w\in F^0$ lies over $v$, then $\phi$ defines a map $F^{\geq 0}_w\to (E_d)^{\geq 0}_v$, which we also denote by $\phi$.

\begin{lemma}\label{lembase}
Let $(F,\phi)$ be an extended representation graph for $E$ and $w$ a vertex in $F$ lying over a vertex $v\in E^0$. Then $(i)-(iii)$ below hold.
\begin{enumerate}[(i)]
\item If $w$ is a source or receives an edge lying over a nonspecial real edge, then $\phi$ defines a bijection $F^{\geq 0}_w\to X_{v}$.
\item If $w$ receives an edge lying over a special real edge $e$, then $\phi$ defines a bijection $F^{\geq 0}_w\to \{x_1\dots x_n\in X_{v}\mid x_1\neq e^*\}$.
\item If $w$ receives an edge lying over a ghost edge, then $\phi$ defines a bijection $F^{\geq 0}_w\to \{x_1\dots x_n\in X_{v}\mid x_1\in (E^1)^*\}$.
\end{enumerate}
In particular, $\phi$ maps paths in $F$ to basis paths in $E_d$. 
\end{lemma}
\begin{proof}
We only prove (i) and leave (ii) and (iii) to the reader. So suppose that $w$ is a source or receives an edge lying over a nonspecial real edge. First we show that $\phi(F^{\geq 0}_w)\subseteq X_v$. Clearly $\phi(w)=v\in X_v$. Now let $p=f_1\dots f_n\in F^{\geq 1}_w$. In order to show that $\phi(p)\in X_v$, it clearly suffices to show that $\phi(f_i)\phi(f_{i+1})\not\in \{ c^*d, e^u(e^u)^*\mid c,d\in E^1, u\in E^0\}$ for any $i\in\{1,\dots,n-1\}$. But that follows from Definition \ref{deferg}(iii),(iv). 

Now we show that the map $F^{\geq 0}_w\to X_{v}$ defined by $\phi$ is injective. We will use the fact that $\phi$ is injective on all sets $s^{-1}(w')~(w'\in F^0)$ (this follows from Definition \ref{deferg}(ii)-(iv)). Suppose that $\phi(p)=\phi(q)$ where $p,q\in F^{\geq 0}_w$. Then $p$ and $q$ have the same length $n$ since $\phi$ preserves lengths. If $n=0$, then $p=w=q$. Now assume that $p=f_1\dots f_n,q=g_1\dots g_n\in F^{\geq 1}_w$. Since $\phi(p)=\phi(q)$, we have $\phi(f_i)=\phi(g_i)~(i=1,\dots,n)$. Clearly $f_1,g_1\in s^{-1}(w)$. Since $\phi$ is injective on $s^{-1}(w)$, it follows that $f_1=g_1$ and hence $f_2,g_2\in s^{-1}(r(f_1))$. Since $\phi$ is injective on $s^{-1}(r(f_1))$, it follows that $f_2=g_2$ and hence $f_3,g_3\in s^{-1}(r(f_2))$. Proceeding like that we obtain $f_i=g_i~(i=1,\dots,n)$ and hence $p=q$ as desired.

It remains to show that the map $F^{\geq 0}_w\to X_{v}$ defined by $\phi$ is surjective. Clearly $X_v=\bigcup_{n\geq 0}X_v^n$ where for any $n\geq 0$, $X_v^n$ is the subset of $X_v$ consisting of all elements of length $n$. Hence it suffices to show that $X_v^n\subseteq \phi(F^{\geq 0}_w)$ for any $n\geq 0$. \\
\\
\underline{$n=0$:} Clearly $X_v^0=\{v\}=\{\phi(w)\}\subseteq\phi(F^{\geq 0}_w)$.\\
\\
\underline{$n=1$:} 
Let $x\in X_v^{1}=s_d^{-1}(v)$. It follows from Definition \ref{deferg}(ii) that $w$ emits an edge $f$ lying over $x$. Hence $X_v^1\subseteq \phi(F^{\geq 0}_w)$.\\
\\
\underline{$n\to n+1$:} 
Assume that $X_v^n\subseteq \phi(F^{\geq 0}_w)$ for some $n\geq 1$. Let $x=x_1\dots x_nx_{n+1}\in X_v^{n+1}$. By the induction assumption we know that there is a $p=f_1\dots f_n\in F^{\geq 0}_w$ such that $\phi(p)=x_1\dots x_n$. It follows from Definition \ref{deferg}(ii)-(iv) that $r(f_n)$ emits an edge $f_{n+1}$ lying over $x_{n+1}$. Hence $\phi(f_1\dots f_nf_{n+1})= x_1\dots x_nx_{n+1}=x$ and thus we have shown that $X_v^{n+1}\subseteq \phi(F^{\geq 0}_w)$.
\end{proof}

\subsubsection{Infinite basis paths}
We denote by $X^{\infty}$ the set of all infinite paths $x=\dots x_3x_2x_1$ in $E_d$ having the property that any finite suffix $x_n\dots x_1$ lies in $X$. We call the elements of $X^{\infty}$ \textit{infinite basis paths}. We call an infinite basis path $x=\dots x_3x_2x_1$ \textit{real} if all of the edges $x_1,x_2, x_3\dots$ are real, and \textit{ghostly} if all of the edges $x_1,x_2, x_3\dots$ are ghost edges. If $x=\dots x_2x_1\in X^{\infty}$ and $n\geq 0$, we set $\tau_{\leq n}(x):=x_{n}\dots x_{1}\in X$ and $\tau_{>n}(x):=\dots x_{n+2}x_{n+1}\in X^{\infty}$. Here we use the convention $\tau_{\leq 0}(x)=s_d(x)$. Two infinite basis paths $x,y\in X^{\infty}$ are called {\it tail-equivalent}, denoted by $x\sim_\infty y$, if there are $m,n\geq 0$ such that $\tau_{>m}(x)=\tau_{>n}(y)$. This defines an equivalence relation on $X^{\infty}$.


\subsubsection{Basis paths $x$ having the property that $x^2$ is again a basis path}
We denote by $X^c$ the set of all finite basis paths $x$ having the property that $x^2$ is again a basis path. Clearly $x\in X^c$ if and only if $x$ is a closed path in $E_d$ that consists either only of real edges or only of ghost edges. We call an $x\in X^c$ \textit{real} if it consists only of real edges, and \textit{ghostly} if it consists only of ghost edges. Two basis paths $x=x_1\dots x_m,y=y_1\dots y_n\in X^c$ are called \textit{equivalent}, denoted by $x\sim_c y$, if $x_1\dots x_m=y_{k+1}\dots y_{n}y_1\dots y_{k}$ for some $1\leq k\leq n$ (note that this implies $m=n$). This defines an equivalence relation on $X^c$. 

\subsection{The connected extended representation graphs containing a source}

Let $v\in E^0$. We define an extended representation graph $(F_v,\phi_v)$ for $E$ by
\begin{align*}
F_v^0=&\{w_x\mid x\in X_v\},\\
F_v^1=&\{f_x\mid x\in X_v\setminus\{v\}\},\\
s_{F_v}(f_x)=&w_{\tau_{\leq |x|-1}(x)},\\
r_{F_v}(f_x)=&w_x,\\
\phi^0_{v}(w_x)=&r_d(x),\\
\phi^1_{v}(f_x)=&\tau_{>|x|-1}(x).
\end{align*}
Note that $F_v$ is connected (since for any $x=x_1\dots x_n\in X_v\setminus\{v\}$ there is a path from $w_v$ to $w_x$, namely $f_{x_1}f_{x_1x_2}\dots f_{x_1\dots x_n}$) and contains a unique source, namely $w_v$.

\begin{example}
Suppose that $E$ is the graph
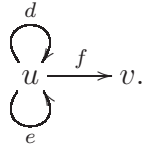
\begin{equation*}
\xymatrix{
 u \ar@(ul,ur)^{d} \ar@(dl,dr)_{e} \ar[r]^{f} &v.
}
\end{equation*}
Then $(F_v,\phi_v)$ is given by 
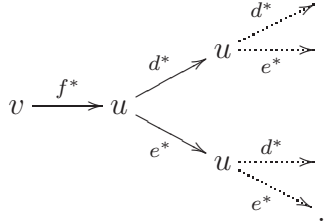
\begin{equation*}
\xymatrix@=10pt{
& & & & &  & & & & &\\
& & & &&u  \ar@{.>}[rr]_{e^*} \ar@{.>}[urr]^{d^*}  \ar@{<-}[dll]_{d^*}&  & &&& \\
& v\ar[rr]^{f^*}& &u & & && & & & \\
& & & &&u\ar@{.>}[rr]^{d^*}   \ar@{.>}[drr]_{e^*}  \ar@{<-}[ull]^{e^*}&&&&&  \\
& & & & &  & & .& & &\\
}
\end{equation*}
Note that in this example $(F_v,\phi_v)$ is a usual representation graph for $E$ as defined in Definition \ref{defrg}. In general, $(F_v,\phi_v)$ is a usual representation graph for $E$ if and only if $v$ is a sink in $E$.
\end{example}

\begin{example}
The extended representation graph in Example \ref{ex1} is isomorphic to $(F_v,\phi_v)$.
\end{example}

\begin{proposition}\label{properg1a}
Let $u,v\in E^0$. Then $(F_u,\phi_u)\cong (F_v,\phi_v)$ if and only if $u=v$.
\end{proposition}
\begin{proof}
Suppose that $\alpha:(F_u,\phi_u)\to (F_v,\phi_v)$ is an isomorphism of extended representation graphs. Let $w$ be the unique source in $F_u$ and $w'$ the unique source in $F_v$. Then $\alpha(w)=w'$ (since a graph homomorphism maps a vertex which is not a source to a vertex which is not a source). It follows that $u=\phi_u(w)=\phi_v(\alpha(w))=\phi_v(w')=v$.
\end{proof}

\begin{lemma}\label{lemerg1}
Let $(F,\phi)$ be an extended representation graph for $E$, $w$ a source in $F$ and $p,q\in F^{\geq 0}_w$. Then $r(p)=r(q)$ if and only $p=q$. 
\end{lemma}
\begin{proof}
Assume that $r(p)=r(q)$. We have to show that this implies $p=q$. First suppose that the length of $p$ or $q$, say $p$, is zero. Then $p=w$. Since $w=r(p)=r(q)$ and $w$ is a source, it follows that the length of $q$ is also zero. Hence $p=q=w$.

Now suppose that $p=f_1\dots f_m,g_1\dots g_n\in F^{\geq 1}_w$. Without loss of generality assume that $m\geq n$. Since $r(p)=r(q)$ and any vertex in $F$ receives at most one edge, we obtain $f_m=g_n,\dots, f_{m-n+1}=g_1$. Assume that $m>n$. Then $f_{m-n}$ is an edge in $F$ ending in $s(f_{m-n+1})=s(g_1)=w$, which contradicts the assumption that $w$ is a source. Hence $m=n$ and $p=q$.
\end{proof}

\begin{proposition}\label{properg1}
Any connected extended representation graph for $E$ that contains a source is isomorphic to one of the extended representation graphs $(F_v,\phi_v)~(v\in E^0)$.
\end{proposition}
\begin{proof}
Suppose that $(F,\phi)$ is a connected extended representation graph for $E$ containing a source $w$ lying over a vertex $v\in E^0$. By Lemma \ref{lembase}(i) there is for any $x\in X_v$ a unique path $p_x$ in $F^{\geq 0}_{w}$ such that $\phi(p_x)=x$. For any $x\in X_v$ set  $w_x:=r(p_x)$. For any $x\in X_v\setminus\{v\}$ let $f_x$ be the last edge of $p_x$. By Lemma \ref{lemerg1} the vertices $w_x~(x\in X_v)$ are pairwise distinct. It follows that the edges  $f_x~(x\in X_v\setminus \{v\})$ are pairwise distinct (since $r(f_x)=r(p_x)=w_x$). 

Let $G$ be the connected subgraph of $F$ defined by $G^0=\{w_x\mid x\in X_v\}$ and $G^1=\{f_x\mid x\in X_v\setminus \{v\}\}$. Then no vertex in $G^0$ receives an edge from $F^1\setminus G^1$ (since $w_v=w$ is a source in $F$ and each $w_x$ where $x\in X_v\setminus \{v\}$ already receives the edge $f_x$). Now let $f\in s_F^{-1}(w_x)$ where $x\in X_v$. Then $p_xf\in F^{\geq 0}_w$. It follows from Lemma \ref{lembase}(i) that $y:=\phi(p_xf)\in \phi(F^{\geq 0}_w)=X_v$. Clearly $p_y=p_xf$ and hence $f=f_y$. Hence no vertex in $G^0$ emits an edge from $F^1\setminus G^1$. It follows that $G$ is a connected component of $F$. This implies that $F=G$ since $F$ is connected. Hence we have
\begin{align*}
F^0=&\{w_x\mid x\in X_v\},\\
F^1=&\{f_x\mid x\in X_v\setminus\{v\}\},\\
s(f_x)=&w_{\tau_{\leq |x|-1}(x)},\\
r(f_x)=&w_x,\\
\phi(w_x)=&r_d(x),\\
\phi(f_x)=&\tau_{>|x|-1}(x)
\end{align*}
and thus $(F,\phi)\cong (F_v,\phi_v)$.
\end{proof}

\subsection{The acyclic and connected extended representation graphs not containing a source}
Let $x=\dots x_3x_2x_1\in X^{\infty}$ denote an infinite basis path. For any $i\in\N$ we denote by $X_i$ the set of all finite basis paths $y=y_1\dots y_n$ of length $\geq 1$ such that $x_iy_1$ is a basis path and $y_1\neq x_{i-1}$ if $i\geq 2$. We define an extended representation graph $(F_{x},\phi_{x})$ for $E$ by 
\begin{align*}
F_{x}^0&=\{w_i\mid i\in \N\}\sqcup \{w_{i,y}\mid i\in \N, y\in X_i\},\\
F_{x}^1&=\{f_i\mid i\in \N\}\sqcup \{f_{i,y}\mid i\in \N, y\in X_i\},\\
s_{F_{x}}(f_i)&=w_{i+1},\\
\quad r_{F_{x}}(f_i)&=w_{i},\\
s_{F_{x}}(f_{i,y})&=\begin{cases}w_{i},\quad&\text{if }|y|=1,\\w_{i,\tau_{\leq |y|-1}(y)},\quad&\text{if }|y|\geq 2,\end{cases}\\
r_{F_{x}}(f_{i,y})&=w_{i,y},\\
\phi_{x}^0(w_i)&=r_d(x_i),\\ 
\phi_{x}^0(w_{i,y})&=r_d(y),\\
\phi_{x}^1(f_i)&=x_i,\\
\phi_{x}^1(f_{i,y})&=\tau_{>|y|-1}(y).
\end{align*}
Note that $F_x$ is connected and does not contain a source. In order to see that $F_x$ is acyclic define a strict partial order $<$ on $F_x^0$ by $w_i<w_j~(i,j\in \N,i>j)$, $w_i<w_{j,y}~(i,j\in\N, i\geq j,y\in X_i)$, and $w_{i,y}<w_{i,z}~(i\in\N, y,z\in X_i, ~y\text{ proper prefix of }z)$. Clearly $s(f)<r(f)$ for any $f\in F_x^1$ and hence $F_x$ is acyclic.

\begin{example}\label{ex3}
Suppose that $E$ is the graph $\xymatrix{ v\ar@[yellow]@(dr,ur)_{e}\ar@[green]@(dl,ul)^{d}}$. Let $x$ be the infinite basis path $\dots d^*d^*d^*$. Then $(F_{x},\phi_{x})$ is given by
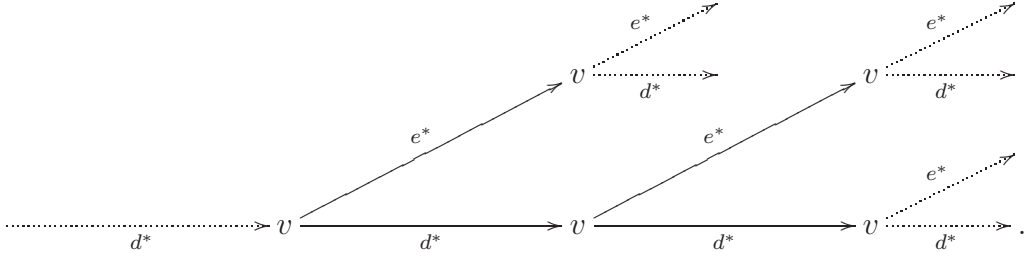
\begin{equation*}
\xymatrix@R=0.7cm@C=0.7cm{
&&&&&&&&&&\ar@[yellow]@{<.}[dll]_{e^*}&&&&\ar@[yellow]@{<.}[dll]_{e^*}\\
&&&&&&&&v\ar@[yellow]@{<-}[ddllll]_{e^*}&&\ar@[green]@{<.}[ll]^{d^*}&&v\ar@[yellow]@{<-}[ddllll]_{e^*}&&\ar@[green]@{<.}[ll]^{d^*}\\
&&&&&&&&&&&&&&\ar@[yellow]@{<.}[dll]_{e^*}\\
\ar@{.>}[rrrr]_{d^*}&&&&v\ar[rrrr]_{d^*}&&&&v&&&&v\ar@[green]@{<-}[llll]^{d^*}&&\ar@[green]@{<.}[ll]^{d^*}.}
\end{equation*}
Note that in this example $(F_{x},\phi_{x})$ is a usual representation graph for $E$. In general, $(F_{x},\phi_{x})$ is a usual representation graph for $E$ if and only if $x$ is ghostly.
\end{example}

\begin{example}\label{ex4}
Suppose that $E$ is the graph $\xymatrix{ v\ar@[yellow]@(dr,ur)_{e}\ar@[green]@(dl,ul)^{d}}$ and choose as special edge $e=e^v$. Let $x$ be the real infinite basis path $\dots ddd$. Then $(F_{x},\phi_{x})$ is given by
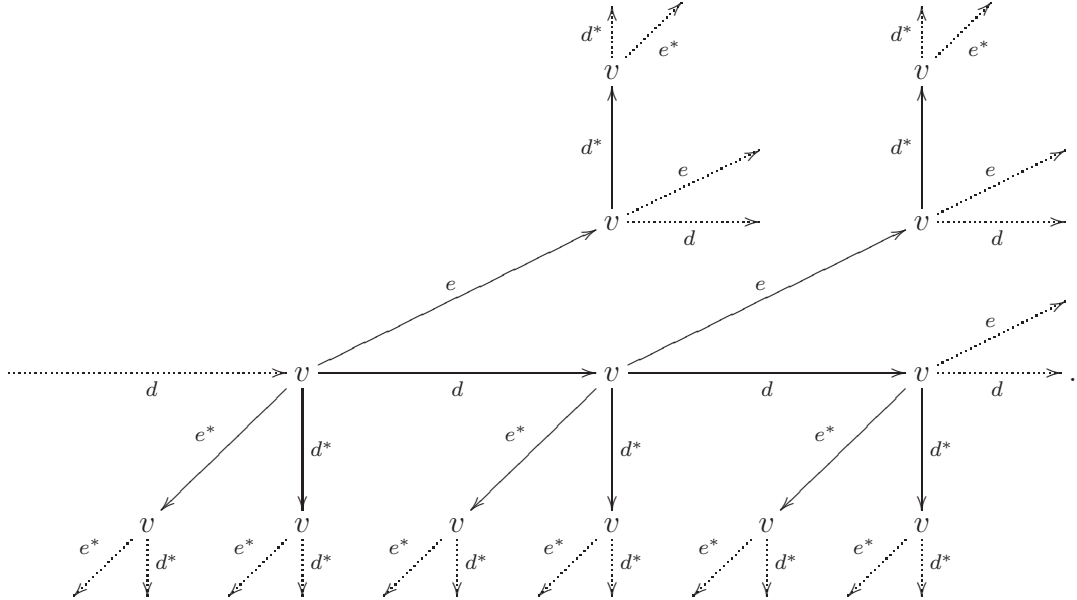
\begin{equation*}
\xymatrix@R=0.7cm@C=0.7cm{
&&&&&&&&&&&&&&\\
&&&&&&&&v\ar@[red]@{.>}[u]^{d^*}\ar@[blue]@{.>}[ur]_{e^*}&&&&v\ar@[red]@{.>}[u]^{d^*}\ar@[blue]@{.>}[ur]_{e^*}&&\\
&&&&&&&&&&\ar@[yellow]@{<.}[dll]_e&&&&\ar@[yellow]@{<.}[dll]_e\\
&&&&&&&&v\ar@[yellow]@{<-}[ddllll]_e\ar@[red][uu]^{d^*}&&\ar@[green]@{<.}[ll]^d&&v\ar@[yellow]@{<-}[ddllll]_e\ar@[red][uu]^{d^*}&&\ar@[green]@{<.}[ll]^d\\
&&&&&&&&&&&&&&\ar@[yellow]@{<.}[dll]_e\\
\ar@{.>}[rrrr]_d&&&&v\ar@[red][dd]^{d^*}\ar@[blue][ddll]_{e^*}\ar[rrrr]_d&&&&v\ar@[red][dd]^{d^*}\ar@[blue][ddll]_{e^*}&&&&v\ar@[green]@{<-}[llll]^d\ar@[red][dd]^{d^*}\ar@[blue][ddll]_{e^*}&&\ar@[green]@{<.}[ll]^d.\\
&&&&&&&&&&&&&&\\
&&v\ar@[red]@{.>}[d]^{d^*}\ar@[blue]@{.>}[dl]_{e^*}&&v\ar@[red]@{.>}[d]^{d^*}\ar@[blue]@{.>}[dl]_{e^*}&&v\ar@[red]@{.>}[d]^{d^*}\ar@[blue]@{.>}[dl]_{e^*}&&v\ar@[red]@{.>}[d]^{d^*}\ar@[blue]@{.>}[dl]_{e^*}&&v\ar@[red]@{.>}[d]^{d^*}\ar@[blue]@{.>}[dl]_{e^*}&&v\ar@[red]@{.>}[d]^{d^*}\ar@[blue]@{.>}[dl]_{e^*}&&\\
&&&&&&&&&&&&&&\\
&&&&&&&&&&&&&&}
\end{equation*}
\end{example}

\begin{proposition}\label{properg2a}
Let $x,z\in X^\infty$. Then $(F_x,\phi_x)\cong (F_z,\phi_z)$ if and only if $x\sim_\infty z$.
\end{proposition}
\begin{proof}
We leave it to the reader to show that $x\sim_\infty z$ implies $(F_x,\phi_x)\cong (F_z,\phi_z)$. In order to show the converse, suppose that $\alpha:(F_x,\phi_x)\to (F_z,\phi_z)$ is an isomorphism of extended representation graphs. We write the vertices and edges in $F_z$ as $w_i'$, $w'_{i,y}$, $f'_i$ and $f'_{i,y}$ to distinguish them from the vertices and edges in $F_x$. Clearly $\alpha$ maps the infinite path $\dots f_3f_2f_1$ in $F_x$ to the infinite path $\dots \alpha(f_3)\alpha(f_2)\alpha(f_1)$ in $F_z$. Because of the strict partial order on $F_z^0$ described right above Example \ref{ex3}, there is an $i\in\N$ such that $r(\alpha(f_i))=w'_j$ for some $j\in\N$. Since for any $k\in\N$, $f'_k$ is the only edge in $F_z$ ending in $w'_k$, it follows that $\dots \alpha(f_{i+2})\alpha(f_{i+1})\alpha(f_{i})=\dots f'_{j+2}f'_{j+1}f'_{j}$. By applying $\phi_z$ to the last equation we obtain $\dots x_{i+2}x_{i+1}x_{i}=\dots z_{j+2}z_{j+1}z_{j}$. Thus $x\sim_\infty z$.
\end{proof}

\begin{lemma}\label{lemerg2}
Let $(F,\phi)$ be an acyclic extended representation graph for $E$, $w\in F^0$ a vertex and $p,q\in F^{\geq 0}_w$. Then $r(p)=r(q)$ if and only $p=q$. 
\end{lemma}
\begin{proof}
Assume that $r(p)=r(q)$. We have to show that this implies $p=q$. First suppose that the length of $p$ or $q$, say $p$, is zero. Then $p=w$. Since $s(q)=w=r(p)=r(q)$ and $F$ is acyclic, it follows that the length of $q$ is also zero. Hence $p=q=w$.

Now suppose that $p=f_1\dots f_m,g_1\dots g_n\in F^{\geq 1}_w$. Without loss of generality assume that $m\geq n$. Since $r(p)=r(q)$ and any vertex in $F$ receives at most one edge, we obtain $f_m=g_n,\dots, f_{m-n+1}=g_1$. Assume that $m>n$. Then $f_1\dots f_{m-n}$ is a path in $F$ starting in $w$ and ending in $s(f_{m-n+1})=s(g_1)=w$, which contradicts the assumption that $F$ is acyclic. Hence $m=n$ and $p=q$.
\end{proof}

\begin{proposition}\label{properg2}
Any nonempty, acyclic and connected extended representation graph for $E$ that does not contain a source is isomorphic to one of the extended representation graphs $(F_{x},\phi_{x})~(x\in X^{\infty})$.
\end{proposition}
\begin{proof}
Suppose that $(F,\phi)$ is a nonempty, acyclic and connected extended representation graph for $E$ that does not contain a source. Choose a vertex $w\in F^0$. Since $F$ does not contain sources, there is a (unique by Condition (i) in Definition \ref{deferg}) infinite path $p=\dots f_3f_2f_1$ in $F$ ending in $w$. Since $F$ is acyclic, the vertices $w_i:=r(f_i)~(i\in \N)$ are pairwise distinct. It follows that the edges $f_i~(i\in\N)$ are pairwise distinct.

Clearly $x:=\phi(p)\in X^\infty$ by Lemma \ref{lembase}. Define the sets $X_i~(i\in\N)$ as in the first paragraph of \S 5.3. Let $i\in\N$. Then, by Lemma \ref{lembase}, there is for any $y\in X_i$ a unique path $p_{i,y}$ in $F^{\geq 1}_{w_i}$ such that $\phi(p_{i,y})=y$. For any $i\in\N$ and $y\in X_i$ set $w_{i,y}:=r(p_{i,y})$ and let $f_{i,y}$ be the last edge of $p_{i,y}$. Suppose that $w_{i,y}=w_{j,z}$ for some $i,j\in\N$, $y\in X_i$ and $z\in X_j$. Assume that $i>j$. Then the paths $p_{i,y}$ and $f_{i-1}\dots f_jp_{j,z}$ have the same source and range. It follows from Lemma \ref{lemerg2} that $p_{i,y}=f_{i-1}\dots f_jp_{j,z}$. By applying $\phi$ to the last equation we obtain $y_1=x_{i-1}$, which contradicts the assumption that $y\in X_i$. We have shown that $i>j$ is impossible. By symmetry $i<j$ is also impossible. It follows that $i=j$. Now Lemma \ref{lemerg2} implies that $y=z$. Thus the vertices $w_{i,y}~(i\in\N,y\in X_i)$ are pairwise distinct. Similarly one can show that the vertices $w_{i,y}~(i\in\N,y\in X_i)$ are distinct from the vertices $w_i~(i\in \N)$. Since $r(f_i)=w_i$ and $r(f_{i,y})=w_{i,y}$, it follows that the edges $f_{i,y}~(i\in\N,y\in X_i)$ are pairwise distinct and distinct from the edges $f_i~(i\in \N)$. 

Let $G$ be the connected subgraph of $F$ defined by $G^0=\{w_i\mid i\in \N\}\sqcup \{w_{i,y}\mid i\in \N, y\in X_i\}$ and $G^1=\{f_i\mid i\in \N\}\sqcup \{f_{i,y}\mid i\in \N, y\in X_i\}$. One checks easily that no vertex in $G^0$ emits or receives an edge from $F^1\setminus G^1$ (cf. the proof of Proposition \ref{properg1}). It follows that $G$ is a connected component of $F$. This implies that $F=G$ since $F$ is connected. Hence we have
\begin{align*}
F^0&=\{w_i\mid i\in \N\}\sqcup \{w_{i,y}\mid i\in \N, y\in X_i\},\\
F^1&=\{f_i\mid i\in \N\}\sqcup \{f_{i,y}\mid i\in \N, y\in X_i\},\\
s(f_i)&=w_{i+1},\\
\quad r(f_i)&=w_{i},\\
s(f_{i,y})&=\begin{cases}w_{i},\quad&\text{if }|y|=1,\\w_{i,\tau_{\leq |y|-1}(y)},\quad&\text{if }|y|\geq 2,\end{cases}\\
r(f_{i,y})&=w_{i,y},\\
\phi(w_i)&=r_d(x_i),\\ 
\phi(w_{i,y})&=r_d(y),
\end{align*}
\begin{align*}
\phi(f_i)&=x_i,\hspace{5.9cm}\\
\phi(f_{i,y})&=\tau_{>|y|-1}(y)
\end{align*}
and thus $(F,\phi)\cong (F_x,\phi_x)$.
\end{proof}

\subsection{The cyclic and connected extended representation graphs}
Let $x=x_1\dots x_m\in X^c$, i.e. $x$ is a closed path in $E_d$ that consists either only of real edges or only of ghost edges. For any $1\leq i\leq m$ we denote by $X_i$ the set of all finite basis paths $y=y_1\dots y_n$ of length $\geq 1$ such that $x_iy_1$ is a basis path and $y_1\neq x_{i+1}$ (here we use the convention $x_{m+1}=x_1$). We define an extended representation graph $(F_x,\phi_x)$ for $E$ by 
\begin{align*}
F_x^0&=\{w_i\mid 1\leq i\leq m\}\sqcup \{w_{i,y}\mid 1\leq i\leq m, y\in X_i\},\\
F_{x}^1&=\{f_i\mid 1\leq i\leq m\}\sqcup \{f_{i,y}\mid 1\leq i\leq m, y\in X_i\},\\
s_{F_{x}}(f_i)&=w_{i-1},\\
\quad r_{F_{x}}(f_i)&=w_{i},\\
s_{F_{x}}(f_{i,y})&=\begin{cases}w_{i},\quad&\text{if }|y|=1,\\w_{i,\tau_{\leq |y|-1}(y)},\quad&\text{if }|y|\geq 2,\end{cases}\\
r_{F_{x}}(f_{i,y})&=w_{i,y},\\
\phi_{x}^0(w_i)&=r_d(x_i),\\ 
\phi_{x}^0(w_{i,y})&=r_d(y),\\
\phi_{x}^1(f_i)&=x_i,\\
\phi_{x}^1(f_{i,y})&=\tau_{>|y|-1}(y)
\end{align*}
using the convention $w_0=w_{m}$. Note that $F_x$ is cyclic (since $c=f_1\dots f_m$ is a cycle), connected and does not contain a source. Define a strict partial order $<$ on $F_x^0$ by $w_i<w_{i,y}~(1\leq i\leq m;~ y\in X_i)$ and $w_{i,y}<w_{i,z}~(1\leq i\leq m; ~y,z\in X_i; ~y\text{ proper prefix of }z)$. Clearly $s(f_{i,y})<r(f_{i,y})$ for any $1\leq i\leq m$ and $ y\in X_i$. Hence the only cycles in $F_x$ are $c=f_1\dots f_m$ and its shifts $f_2\dots f_m f_1$, $f_3\dots f_mf_1f_2$ etc. 

\begin{example}\label{ex1.75}
Suppose that $E$ is the graph $\xymatrix{ v\ar@[yellow]@(dr,ur)_{e}\ar@[green]@(dl,ul)^{d}}$ and let $x=e^*$. Then $(F_{x},\phi_{x})$ is given by
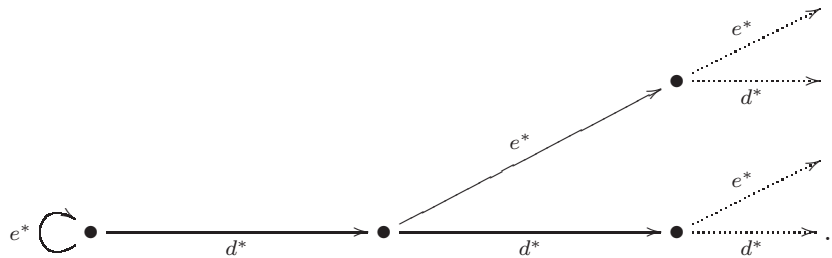
\begin{equation*}
\xymatrix@R=0.7cm@C=0.7cm{
&&&&&&&&&&&\ar@[yellow]@{<.}[dll]_{e^*}\\
&&&&&&&&&\bullet\ar@[yellow]@{<-}[ddllll]_{e^*}&&\ar@[green]@{<.}[ll]^{d^*}\\
&&&&&&&&&&&\ar@[yellow]@{<.}[dll]_{e^*}\\
&\bullet\ar@[yellow]@(dl,ul)^{e^*}&&&&\bullet\ar@[green]@{<-}[llll]^{d^*}&&&&\bullet\ar@[green]@{<-}[llll]^{d^*}&&\ar@[green]@{<.}[ll]^{d^*}.}
\end{equation*}
Note that in this example $(F_{x},\phi_{x})$ is a usual representation graph for $E$. In general, $(F_{x},\phi_{x})$ is a usual representation graph for $E$ if and only if $x$ is ghostly.
\end{example}

\begin{example}\label{ex2}
Suppose that $E$ is the graph $\xymatrix{ v\ar@[yellow]@(dr,ur)_{e}\ar@[green]@(dl,ul)^{d}}$, choose $e=e^v$ as special edge and let $x=e$. Then $(F_{x},\phi_{x})$ is given by
\newpage
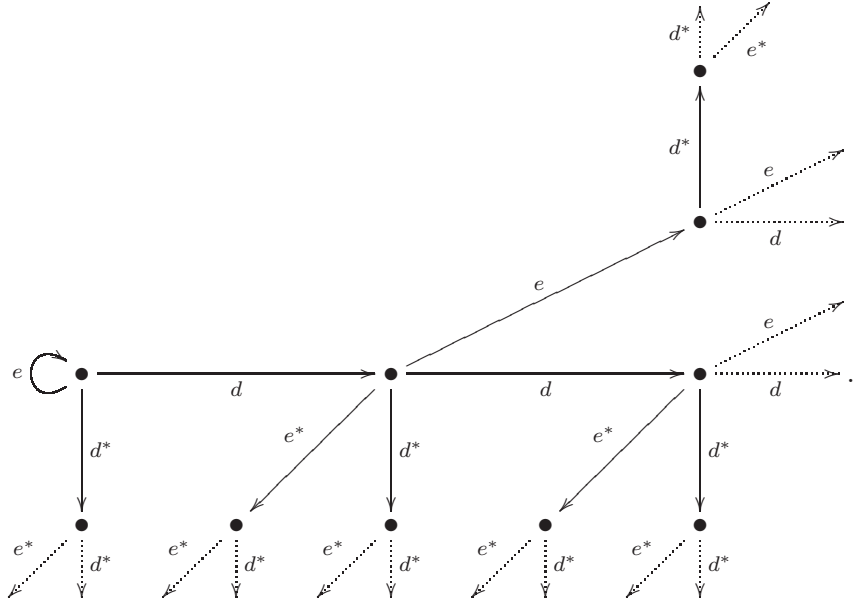
\begin{equation*}
\xymatrix@R=0.7cm@C=0.7cm{
&&&&&&&&&&&\\
&&&&&&&&&\bullet\ar@[red]@{.>}[u]^{d^*}\ar@[blue]@{.>}[ur]_{e^*}&&\\
&&&&&&&&&&&\ar@[yellow]@{<.}[dll]_e\\
&&&&&&&&&\bullet\ar@[yellow]@{<-}[ddllll]_e\ar@[red][uu]^{d^*}&&\ar@[green]@{<.}[ll]^d\\
&&&&&&&&&&&\ar@[yellow]@{<.}[dll]_e\\
&\bullet\ar@[yellow]@(dl,ul)^e\ar@[red][dd]^{d^*}&&&&\bullet\ar@[green]@{<-}[llll]^d\ar@[red][dd]^{d^*}\ar@[blue][ddll]_{e^*}&&&&\bullet\ar@[green]@{<-}[llll]^d\ar@[red][dd]^{d^*}\ar@[blue][ddll]_{e^*}&&\ar@[green]@{<.}[ll]^d.\\
&&&&&&&&&&&\\
&\bullet\ar@[red]@{.>}[d]^{d^*}\ar@[blue]@{.>}[dl]_{e^*}&&\bullet\ar@[red]@{.>}[d]^{d^*}\ar@[blue]@{.>}[dl]_{e^*}&&\bullet\ar@[red]@{.>}[d]^{d^*}\ar@[blue]@{.>}[dl]_{e^*}&&\bullet\ar@[red]@{.>}[d]^{d^*}\ar@[blue]@{.>}[dl]_{e^*}&&\bullet\ar@[red]@{.>}[d]^{d^*}\ar@[blue]@{.>}[dl]_{e^*}&&\\
&&&&&&&&&&&\\
&&&&&&&&&&&}
\end{equation*}
\end{example}

\begin{proposition}\label{properg3a}
Let $x,z\in X^c$. Then $(F_x,\phi_x)\cong (F_z,\phi_z)$ if and only if $x\sim_c z$.
\end{proposition}
\begin{proof}
We leave it to the reader to show that $x\sim_c z$ implies $(F_x,\phi_x)\cong (F_z,\phi_z)$. In order to show the converse, suppose that $\alpha:(F_x,\phi_x)\to (F_z,\phi_z)$ is an isomorphism of extended representation graphs. Write $x=x_1\dots x_m$ and $z=z_1\dots z_n$. We write the edges in $F_z$ as $f'_i$ and $f'_{i,y}$ to distinguish them from the edges in $F_x$. Clearly $\alpha$ maps the cycle $f_1\dots f_m$ in $F_x$ to the cycle $\alpha(f_1)\dots \alpha(f_m)$ in $F_z$. By the last sentence before Example \ref{ex1.75} there is a $1\leq k\leq n$ such that $\alpha(f_1)\dots \alpha(f_m)=f'_{k+1}\dots f'_nf'_1\dots f'_k$. By applying $\phi_z$ to the last equation we obtain $x_1\dots x_m=z_{k+1}\dots z_nz_1\dots z_k$. Thus $x\sim_c z$.
\end{proof}

\begin{lemma}\label{lemerg3}
Let $(F,\phi)$ be an extended representation graph for $E$, $w\in F^0$ a vertex and $p,q\in F^{\geq 1}_w$. If $r(p)=r(q)$, then either $p=q$, $cp=q$ or $p=cq$ where $c$ is some closed path in $F$. 
\end{lemma}
\begin{proof}
Assume that $r(p)=r(q)$. Write $p=f_1\dots f_m$ and $q=g_1\dots g_n\in F^{\geq 1}_w$. Assume that $m\geq n$. Since $r(p)=r(q)$ and any vertex in $F$ receives at most one edge, we obtain $f_m=g_n,\dots, f_{m-n+1}=g_1$. If $m=n$, then $p=q$. Assume now that $m>n$. Then $c:=f_1\dots f_{m-n}$ is a path in $F$ starting in $w$ and ending in $s(f_{m-n+1})=s(g_1)=w$. Clearly $p=cq$. Similarly one can show that $cp=q$ for some closed path $c$ if $m<n$.
\end{proof}

\begin{proposition}\label{properg3}
Any cyclic and connected extended representation graph for $E$ is isomorphic to one of the extended representation graphs $(F_{x},\phi_{x})~(x\in X^{c})$.
\end{proposition}
\begin{proof}
Suppose that $(F,\phi)$ is a cyclic and connected extended representation graph for $E$. Choose a cycle $c=f_1\dots f_m$ in $F$. Since $c$ is a cycle, the vertices $w_i:=r(f_i)~(1\leq i \leq m)$ are pairwise distinct. It follows that the edges $f_i~(1\leq i \leq m)$ are pairwise distinct.

Clearly $x:=\phi(c)\in X^c$ by Lemma \ref{lembase}. Define the sets $X_i~(1\leq i \leq m)$ as in the first paragraph of \S 5.4. Let $1\leq i \leq m$. Then, by Lemma \ref{lembase}, there is for any $y\in X_i$ a unique path $p_{i,y}$ in $F^{\geq 1}_{w_i}$ such that $\phi(p_{i,y})=y$. For any $1\leq i \leq m$ and $y\in X_i$ set $w_{i,y}:=r(p_{i,y})$ and let $f_{i,y}$ be the last edge of $p_{i,y}$. Suppose that $w_{i,y}=w_{j,z}$ for some $1\leq i,j \leq m$, $y\in X_i$ and $z\in X_j$. Assume that $i<j$. Then the paths $p_{i,y}$ and $f_{i+1}\dots f_jp_{j,z}$ have the same source and range. It follows from Lemma \ref{lemerg3} that 
\begin{align}
p_{i,y}&=f_{i+1}\dots f_jp_{j,z} \text{ or}\\
cp_{i,y}&=f_{i+1}\dots f_jp_{j,z} \text{ or}\\
p_{i,y}&=cf_{i+1}\dots f_jp_{j,z}
\end{align}
where $c$ is a closed path at $w_i$. Clearly $c$ must be a power of the cycle $f_{i+1}\dots f_i$. If (6) or (8) holds, then we obtain $y_1=x_{i+1}$ by applying $\phi$ to the corresponding equation. But that contradicts the assumption that $y\in X_i$. If (7) holds, then we obtain $x_{j+1}=z_1$ by applying $\phi$. But that contradicts the assumption that $z\in X_j$. Hence $i<j$ is impossible. By symmetry $i>j$ is also impossible. It follows that $i=j$. Hence the paths $p_{i,y}$ and $p_{j,z}$ have the same source and range. By Lemma \ref{lemerg3} we have 
\begin{align}
p_{i,y}&=p_{j,z} \text{ or}\\
cp_{i,y}&=p_{j,z} \text{ or}\\
p_{i,y}&=cp_{j,z}
\end{align}
where $c$ is a closed path at $w_i$. Clearly $c$ must be a power of the cycle $f_{i+1}\dots f_i$. If (10) holds, then we obtain $x_{i+1}=z_1$ by applying $\phi$. But that contradicts the assumption that $z\in X_j=X_i$. If (11) holds, then we obtain $y_1=x_{i+1}$. But that contradicts the assumption that $y\in X_i$. Hence (9) holds. By applying $\phi$ to (9) we obtain $y=z$ as desired. We have shown that the vertices $w_{i,y}~(1\leq i \leq m,y\in X_i)$ are pairwise distinct. Similarly one can show that the vertices $w_{i,y}~(1\leq i \leq m,y\in X_i)$ are distinct from the vertices $w_i~(1\leq i \leq m)$. Since $r(f_i)=w_i$ and $r(f_{i,y})=w_{i,y}$, it follows that the edges $f_{i,y}~(1\leq i \leq m,y\in X_i)$ are pairwise distinct and distinct from the edges $f_i~(1\leq i \leq m)$. 

Let $G$ be the connected subgraph of $F$ defined by $G^0=\{w_i\mid 1\leq i \leq m\}\sqcup \{w_{i,y}\mid 1\leq i \leq m, y\in X_i\}$ and $G^1=\{f_i\mid 1\leq i \leq m\}\sqcup \{f_{i,y}\mid 1\leq i \leq m, y\in X_i\}$. One checks easily that no vertex in $G^0$ emits or receives an edge from $F^1\setminus G^1$ (cf. the proof of Proposition \ref{properg1}). It follows that $G$ is a connected component of $F$. This implies that $F=G$ since $F$ is connected. Hence we have
\begin{align*}
F^0&=\{w_i\mid 1\leq i\leq m\}\sqcup \{w_{i,y}\mid 1\leq i\leq m, y\in X_i\},\\
F^1&=\{f_i\mid 1\leq i\leq m\}\sqcup \{f_{i,y}\mid 1\leq i\leq m, y\in X_i\},\\
s(f_i)&=w_{i-1},\\
\quad r(f_i)&=w_{i},\\
s(f_{i,y})&=\begin{cases}w_{i},\quad&\text{if }|y|=1,\\w_{i,\tau_{\leq |y|-1}(y)},\quad&\text{if }|y|\geq 2,\end{cases}\\
r(f_{i,y})&=w_{i,y},\\
\phi(w_i)&=r_d(x_i),\\ 
\phi(w_{i,y})&=r_d(y),\\
\phi(f_i)&=x_i,\\
\phi(f_{i,y})&=\tau_{>|y|-1}(y)
\end{align*}
and thus $(F,\phi)\cong (F_x,\phi_x)$.
\end{proof}

\subsection{Summary}

The theorem below follows from Propositions \ref{properg1a}, \ref{properg1}, \ref{properg2a}, \ref{properg2}, \ref{properg3a} and \ref{properg3}.

\begin{theorem}\label{thmm1}
Let $R$ (resp. $S$) be a complete set of representatives for the $\sim_\infty$-equivalence classes (resp. $\sim_c$-equivalence classes). Then $\{(F_{v},\phi_{v}),(F_{x},\phi_{x}),(F_{y},\phi_{y})\mid v\in E^0, x\in R, y\in S\}$ is a complete set of representatives for the isomorphism classes of the nonempty and connected extended representation graphs for $E$.
\end{theorem}

We call a $\sim_\infty$-equivalence class or $\sim_c$-equivalence class \textit{ghostly} if it contains a ghostly element (note that in this case any element of this class is ghostly). In view of Example \ref{extrivial} we obtain the following classification of the usual representation graphs (which, in turn, yields a classification of the usual algebraic branching systems in view of \S 3.3).

\begin{theorem}\label{corm1}
Let $R_{\ghost}$ (resp. $S_{\ghost}$) be a complete set of representatives for the ghostly $\sim_\infty$-equivalence classes (resp. ghostly $\sim_c$-equivalence classes). Then $\{(F_{v},\phi_{v}),(F_{x},\phi_{x}),(F_{y},\phi_{y})\mid v\in E^0_{\sink}, x\in R_{\ghost}, y\in S_{\ghost}\}$ is a complete set of representatives for the isomorphism classes of the nonempty and connected representation graphs for $E$.
\end{theorem}

\section{A class of counter-examples to the converse of Schur's lemma}

\begin{theorem}\label{thmschur}
Let $x$ be a cycle in $E$ which consists only of special edges and has an exit. Then the $L(E)$-module $W:=W(F_x,\phi_x)$ is not simple and $\End_{L(E)}(W)\cong K$.
\end{theorem}
\begin{proof}
Write $x=x_1\dots x_m$ and choose an exit $e$ of $x$. Then $s(e)=s(x_i)$ and $e\neq x_i$ for some $1\leq i\leq m$. Since $x_i$ is special, $e$ is not special. Clearly $e\in X_{i-1}$ (recall that $X_{i-1}$ is the set of all basis paths $y=y_1\dots y_n$ of length $\geq 1$ such that $x_{i-1}y_1$ is a basis path and $y_1\neq x_{i}$). Let $V$ be the $K$-span of the vertices $w_{i-1,y}~(y\in X_{i-1}, y_1=e)$. One checks easily that $V$ is a submodule of $W$. Since $V$ contains $w_{i-1,e}$ but not $w_{i-1}$, it is a proper submodule of $W$. Hence $W$ is not simple.

It remains to show that $\End_{L(E)}(W)\cong K$. Let $\theta\in \End_{L(E)}(W)$. Then 
\begin{equation}
\theta(w_m)=\sum_{1\leq i\leq m}k_iw_i+\sum_{1\leq i\leq m,y\in X_i}k_{i,y}w_{i,y}
\end{equation}
for some $k_i,k_{i,y}\in K~(1\leq i\leq m, y\in X_{i})$ of which only a finite number are nonzero. Clearly $\theta(w_m)=\theta(w_m.x)=\theta(w_m).x$. Hence we obtain
\begin{equation}
\sum_{1\leq i\leq m}k_iw_i+\sum_{1\leq i\leq m,y\in X_i}k_{i,y}w_{i,y}=\sum_{1\leq i\leq m}k_iw_i.x+\sum_{1\leq i\leq m,y\in X_i}k_{i,y}w_{i,y}.x.
\end{equation}
Clearly 
\begin{equation}
w_i.x=\delta_{i,m}w_i
\end{equation}
for any $1\leq i\leq m$. On the other hand, we have for any $1\leq i\leq m$ and $y\in X_i$\newpage 
\begin{equation}
w_{i,y}.x=\begin{cases}
w_{i,px},\quad &\text{if }y=p\text{ and }r(p)=s(x),\\
w_{i,q^*},\quad &\text{if }y=q^*x^*,\\
w_{i,pq^*},\quad &\text{if }y=pq^*x^*,\\
w_{i,px_{k+1}\dots x_m},\quad &\text{if }y=px_k^*\dots x_1^*\text{ for some }1\leq k \leq m,\\
0,\quad &\text{otherwise},
\end{cases}
\end{equation}
where $p,q\in E^{\geq 1}$. Note that $y_1\neq x_k^*$ for any $1\leq k\leq m$, since $x_iy_1$ is a basis path, $x$ is a cycle and $x_i$ is special. It follows from (13), (14) and (15) that 
\begin{equation}
k_i=0 ~(1\leq i\leq m-1)
\end{equation}
and 
\begin{equation}
\sum_{y\in X_i}k_{i,y}w_{i,y}=\sum_{y\in X_i}k_{i,y}w_{i,y}.x~(1\leq i\leq m).
\end{equation}
Fix a $1\leq i \leq m$ and set $Y_i:=\{y\in X_i\mid k_{i,y}\neq 0\}$. Note that $Y_i$ is a finite set. It follows from equation (17) that
\begin{equation}
\sum_{y\in Y_i}k_{i,y}w_{i,y}=\sum_{y\in Y_i}k_{i,y}w_{i,y}.x.
\end{equation}
Clearly $w_{i,y}.x\neq 0$ for any $y\in Y_i$, otherwise the right hand side of (18) would have less summands than the left hand side. Similarly $w_{i,y}.x\neq w_{i,z}.x$ for any $y\neq z\in Y_i$. For any $n\geq 0$ we denote by $Y_i^n$ the subset of $Y_i$ consisting of all $y$ which have precisely $n$ letters that are ghost edges (of course these must be the last $n$ letters). We show by induction that $Y_i^n=\emptyset$ for any $n\geq 0$.\\
\\
\underline{$n=0$}: Assume that $Y^0_i$ is not the empty set. Then we can choose an element $y\in Y^0_i$ with maximal length. Clearly $y=p$ for some $p\in E^{\geq 1}$ and hence $w_{i,y}.x=w_{i,px}$ by (15). Therefore $w_{i,px}$ has to appear with nonzero coefficient also on the left hand side of (18), i.e. $k_{i,px}\neq 0$. It follows that $px\in Y^0_i$, which is impossible since $y=p$ has maximal length among the elements of $Y^0_i$. Thus $Y^0_i=\emptyset$.\\
\\
\underline{$n\to n+1$}: Assume that $Y^{0}_i=\dots=Y^n_i=\emptyset$ for some $n\geq 0$. Suppose that $Y^{n+1}_i\neq\emptyset$. Choose a $y\in Y^{n+1}_i$. By (15), $w_{i,y}.x=w_{i,z}$ for some $z\in X_i$. Clearly $z$ has less letters which are ghost edges than $y$. By the induction assumption it follows that $z\not\in Y_i$ and hence $k_{i,z}=0$. But since $w_{i,y}.x=w_{i,z}$, the vertex $w_{i,z}$ has to appear with nonzero coefficient on the left hand side of (18), a contradiction. Thus $Y^{n+1}_i=\emptyset$.

We have shown that $Y_i^n=\emptyset$ for any $n\geq 0$. Since $Y_i=\bigcup_{n\geq 0}Y_i^n$, it follows that $Y_i=\emptyset$. Hence 
\begin{equation}
k_{i,y}=0~(1\leq i\leq m,y\in X_i).
\end{equation} 
It follows from (12), (16) and (19) that $\theta(w_m)=k_m w_m$. Hence for any $\theta\in\End_{L(E)}(W)$ there is a $k(\theta)\in K$ such that $\theta(w_m)=k(\theta) w_m$. Let $\xi:\End_{L(E)}(W)\to K$ be the map defined by $\xi(\theta)=k(\theta)$. One checks routinely that $\xi$ is a ring isomorphism.
\end{proof}

\begin{remark}
The following example shows that the assumption in Theorem \ref{thmschur} that the cycle $x$ has an exit is necessary. Suppose that $E$ is the graph $\xymatrix{v\ar@[yellow]@(dr,ur)_{e}}$ and $x=e$. Then $(F_x,\phi_x)$ is given by $\xymatrix{v\ar@[yellow]@(dr,ur)_{e}}$. Clearly the module $W(F_x,\phi_x)$ is simple.
\end{remark}

Suppose $x$ and $x'$ are cycles in $E$ which consist only of special edges and have an exit. If $x\sim_c x'$, then $(F_x,\phi_x)\cong(F_{x'},\phi_{x'})$ by Proposition \ref{properg3a}. It follows that the modules $W(F_x,\phi_x)$ and $W(F_{x'},\phi_{x'})$ are isomorphic since $W$ is a functor. One can ask if $W(F_x,\phi_x)$ and $W(F_{x'},\phi_{x'})$ can be isomorphic if $x\not\sim_c x'$. The theorem below gives a partial answer to this question.

\begin{theorem}\label{thmschur2}
Let $x$ and $x'$ be cycles in $E$ which do not have a common vertex. Then the $L(E)$-modules $W(F_x,\phi_x)$ and $W(F_{x'},\phi_{x'})$ are not isomorphic.
\end{theorem}
\begin{proof}
Suppose that there is an isomorphism $\theta:W(F_{x'},\phi_{x'})\to W(F_x,\phi_x)$. Write $x=x_1\dots x_m$ and $x'=x'_1\dots x'_{m'}$. Moreover, write the vertices in $F_{x'}$ as $w'_i$ and $w'_{i,y}$ to distinguish them from the vertices in $F_x$. Clearly 
\begin{equation}
\theta(w'_{m'})=\sum_{1\leq i\leq m}k_iw_i+\sum_{1\leq i\leq m,y\in X_i}k_{i,y}w_{i,y}
\end{equation}
for some $k_i,k_{i,y}\in K~(1\leq i\leq m, y\in X_{i})$ of which only a finite number are nonzero. Since $\theta(w'_{m'})=\theta(w'_{m'}.x')=\theta(w'_{m'}).x'$, it follows that
\begin{equation}
\sum_{1\leq i\leq m}k_iw_i+\sum_{1\leq i\leq m,y\in X_i}k_{i,y}w_{i,y}=\sum_{1\leq i\leq m}k_iw_i.x'+\sum_{1\leq i\leq m,y\in X_i}k_{i,y}w_{i,y}.x'.
\end{equation}
Clearly 
\begin{equation}
w_i.x'=0
\end{equation}
for any $1\leq i\leq m$ since $x$ and $x'$ have no common vertex. On the other hand, we have for any $1\leq i\leq m$ and $y\in X_i$
\begin{equation}
w_{i,y}.x'=\begin{cases}
w_{i,px'},\quad &\text{if }y=p\text{ and }r(p)=s(x'),\\
w_{i,q^*},\quad &\text{if }y=q^*(x')^*,\\
w_{i,pq^*},\quad &\text{if }y=pq^*(x')^*,\\
w_{i,px'_{k+1}\dots x'_{m'}},\quad &\text{if }y=p(x'_k)^*\dots (x'_1)^*\text{ for some }1\leq k \leq m',\\
0,\quad &\text{otherwise},
\end{cases}
\end{equation}
where $p,q\in E^{\geq 1}$. Note that $y_1\neq (x'_k)^*$ for any $1\leq k\leq m'$, since $x_iy_1$ is a basis path and $x$ and $x'$ do not have a common vertex. It follows from (21), (22) and (23) that
\begin{equation}
k_i=0 ~(1\leq i\leq m)
\end{equation}
and 
\begin{equation}
\sum_{y\in X_i}k_{i,y}w_{i,y}=\sum_{y\in X_i}k_{i,y}w_{i,y}.x'~(1\leq i\leq m).
\end{equation}
One can deduce from (25) that
\begin{equation}
k_{i,y}=0~(1\leq i\leq m,y\in X_i).
\end{equation} 
The proof is essentially the same as the proof that (17) implies (19) (just replace $x$ by $x'$). It follows from (20), (24) and (26) that $\theta(w'_{m'})=0$, a contradiction. Thus $W(F_x,\phi_x)$ and $W(F_{x'},\phi_{x})$ are not isomorphic.
\end{proof}

\section*{Acknowledgments}
I would like to thank Daniel Gon\c{c}alves for encouraging me to write this paper. I would also like to thank Tran Giang Nam, who asked me if I know any modules over Leavitt path algebras that are counter-examples to the converse of Schur's lemma. This question led to Section 6.

\bibliographystyle{alpha} 
\bibliography{Bibb}

\end{document}